
\input amstex
\input amsppt.sty
\NoBlackBoxes
\leftheadtext{K. Matsuki and R. Wentworth}
\rightheadtext{Vector bundles on an algebraic surface}
\topmatter
\title Mumford-Thaddeus principle on the moduli
space\linebreak
of vector bundles on an algebraic surface \endtitle
\author Kenji Matsuki and Richard Wentworth\endauthor
\endtopmatter

\document

$$\bold{\text{Contents}}$$

\S 0. Introduction

\S 1. Some Basic Finiteness Results

\S 2. Key GIT Lemma after Simpson

\S 3. Rationally Twisted Gieseker-Semistability and Stratifications

\S 4. Construction of Flip: Integral Case

\S 5. Construction of Flip: Rational Case

\S 6. References

\vskip.3in

\S 0. Introduction.

\vskip.1in

The purpose of this paper is to study what we call the
``Mumford-Thaddeus principle" which states that Geometric Invariant Theory
 (henceforth, ``GIT") quotients undergo specific transformations (birational
and
similar to Mori's flip under some mild conditions) when the polarization
 (i.e. the linearized ample line bundle) is varied (cf.[MFK94],
[Thaddeus93,94],
 and [Dolgachev-Hu93]).  The case we consider is the moduli space of vector
 bundles on an algebraic surface.

Let $X$ be a nonsingular projective surface over an
algebraically closed field $k$ of characteristic zero
(this assumption about characteristic is used only in analyzing the behavior of
semistability with respect to Galois covers in \S 5).  To construct
 the moduli space, [Gieseker77] first introduces the notion of
$H$-Gieseker-semistability for torsion-free coherent sheaves; a refinement
 of Mumford-Takemoto $H$-slope-semistability.  This depends on the
 choice of an ample line bundle $H$ (actually, only the numerical class
 of $H$) on the surface $X$.  Then he considers an appropriate GIT problem
 and proves that the projective quotient $M(r,c_1,c_2,H)$ coarsely represents
the
Seshadri equivalence classes of torsion-free coherent sheaves of
 fixed rank $r$, first and second Chern classes $c_1$ and $c_2$, which are
$H$-Gieseker-semistable.  Once the existence of such a moduli space
is established, a basic and natural problem is to study the change of
 $M(r,c_1,c_2,H)$ as $H$ varies.

\vskip.1in

The aim of this paper is to prove the

\proclaim{Main Theorem} As the polarization $H$ varies through the ample
 cone of $X$, the moduli space $M(r,c_1,c_2,H)$ goes through a locally finite
 sequence of Thaddeus-type flips (and contraction morphisms and their inverses)
in the category of projective
 moduli spaces $M((r,c_1,c_2) \otimes
{\Cal L},A)$ of ${\Cal L}$-twisted
$A$-Gieseker-semistable sheaves for rational line bundles
${\Cal L} \in \text{Pic}(X) \otimes {\Bbb Q}$ and ample line bundles $A$
 (see below
for the definition of ${\Cal L}$-twisted semistability and the definition of a
Thaddeus-type flip).
\endproclaim

In this paper a Thaddeus-type flip is defined to be a diagram
$$
\CD
{Q_H}^{ss}//G @.@.@.@. {Q_{H'}}^{ss}//G\\
@.\psi \searrow @.@.\swarrow \psi^+@.\\
@.@.{Q_A}^{ss}//G@.@.\\
\endCD
$$
where ${Q_H}^{ss}//G$ (resp. ${Q_{H'}}^{ss}//G$, ${Q_A}^{ss}//G$) is the GIT
quotient
of the locus of semistable points ${Q_H}^{ss}$ (resp. ${Q_{H'}}^{ss}$,
${Q_A}^{ss}$) with
respect to the linearization $H$ (resp. $H'$, $A$) of the action of a reductive
group $G$
 on a projective scheme $Q$, and the morphisms $\psi$ and $\psi^+$ are induced
from the inclusion
 ${Q_H}^{ss} \subset {Q_A}^{ss} \supset {Q_{H'}}^{ss}$.  We remark that we do
not require the morphisms
 to be birational or small in the definition, though [Dolgachev-Hu93] and
[Thaddeus94] claim these and other nice properties for the morphisms under some
mild conditions.

It should be mentioned that there exists another compactification of the moduli
space of vector bundles,
 the Uhlenbeck compactification, whose projective
 variety structure is described in [J.Li93].  Its variation with respect to
 the polarization has been considered, at least in rank 2 case, in
[Hu-W.P.Li94].
  The introduction in the present paper of the notion of
${\Cal L}$-twisted $A$-Gieseker-semistable sheaves is the novelty required to
study
 Gieseker's space.

\vskip.1in

We now outline the structure of this paper:
In \S 1 we prove some basic finiteness results.
[Gieseker77] shows that if we fix the polarization $H$, then the
set $S(r,c_1,c_2,H)$ of all coherent torsion-free sheaves
with given numerical data and which are $H$-Gieseker-semistable is bounded,
i.e., parametrized by a scheme of finite type
over $k$.  Global finiteness, by contrast, does not hold in general,
 i.e., the union of the sets $S(r,c_1,c_2,H)$ when $H$ varies over all ample
 numerical classes may not be bounded (see Remark 1.7 (1) for an example).
The best one can hope for is local finiteness, and we prove this for slope
semistability
 (which implies the result for Gieseker-semistability): Let $\Delta$ be a
convex polyhedral
 cone in the ample cone of $X$ and $\mu(r,c_1,c_2,H)$ denote the set of all
$H$-slope-semistable torsion-free sheaves of given numerical data.  Then the
union of the sets
$\mu(r,c_1,c_2,H)$ when $H$ varies over $\Delta - \{0\}$ is bounded.

We also have the following finiteness result: The set of all subsheaves $F
\subset
E$ where $E \in \mu(r,c_1,c_2,H)$
 for some $H \in \Delta - \{0\}$ with torsion-free quotient $E/F$ satisfying
$$\{\frac{c_1(F)}{rk(F)} - \frac{c_1(E)}{rk(E)}\} \cdot H = 0$$
is bounded.

Hence the set of all hyperplanes
$$L_F := \{z \in N^1(X)_{\Bbb Q};\{\frac{c_1(F)}{rk(F)} -
\frac{c_1(E)}{rk(E)}\} \cdot z = 0\}$$
is finite, where $F$ runs over all such subsheaves as above with
$$\frac{c_1(F)}{rk(F)} - \frac{c_1(E)}{rk(E)}$$
being not numerically trivial.  The $L_F$'s thus give a stratification $\Delta
-
\{0\} = \coprod_s\Delta_s$ into cells
 which characterizes the change of the set $\mu(r,c_1,c_2,H)$.

\vskip.1in

In \S 2, we study the GIT problem.  Ideally, one would
like to use the boundedness result above to embed all the
sheaves in $\cup_{H \in \Delta - \{0\}}S(r,c_1,c_2,H)$
 into one big space independent of $H$ (e.g. a
Grothendieck Quot scheme), observe the correspondence
between $H$-Gieseker-semistable sheaves and GIT semistable
points for an appropriate linearization, and then apply the
Mumford-Thaddeus principle.  Though this turns out to be
too naive a picture and indeed false if one applies the
known construction directly,  it will nevertheless be our
guiding principle throughout.

Gieseker's original construction using the covariant
method is not suitable for our purpose since the relevant
spaces are highly dependent of the polarization, whereas
Simpson's direct use of the Grothendieck Quot scheme
[Simpson92] provides us with a more adequate setting.
Briefly, the method for a FIXED polarization is the
following: First, one embeds  all the sheaves $E \in
S(r,c_1,c_2,H)$ into the Quot scheme parametrizing
surjections ${\Cal O}_X^{\oplus l} \rightarrow E \otimes
H^a \rightarrow 0$ for large $a$.  Second, one proves the
correspondence between $H$-Gieseker-semistable sheaves and
the GIT semistable points  with respect to the action of
$SL(l)$ and the linearization induced from the embedding
of the Quot scheme into a Grassmannian with its Pl\"ucker
coordinates.  The latter embedding comes from taking a
further high multiple of $H$ and a surjection $H^0({\Cal
O}_X^{\oplus l} \otimes H^m) \rightarrow H^0(E \otimes H^a
\otimes H^m) \rightarrow 0$.  The GIT quotient then
recovers $M(r,c_1,c_2,H)$.

Our Key GIT Lemma (Lemma 2.4), which is a modification of
Simpson's method, goes along the same line -- one simply
uses DIFFERENT ample bundles $A$ and $H$ for the first and
second embeddings -- but this turns out to be  quite
subtle and provides us with some unexpected consequences.
First, by the boundedness result it is a simple matter to
embed all the sheaves in $\cup_{H \in \Delta -
\{0\}}S(r,c_1,c_2,H)$ into a Quot scheme by taking a high
multiple of $A \in \Delta - \{0\}$ and a surjection ${\Cal
O}_X^{\oplus l} \rightarrow E \otimes A^a \rightarrow
0$, though it will be observed that this only reflects the
properties local around $A$.  If one then considers the
linearization obtained by embedding the Quot scheme into
a Grassmannian via a surjection $H^0({\Cal O}_X^{\oplus
l} \otimes H^m) \rightarrow H^0(E \otimes A^a \otimes
H^m) \rightarrow 0$, the question of what sheaves the
semistable points correspond to becomes quite delicate.
The answer is the following: the set of semistable
points with respect to the second embedding by $H$
corresponds not to $S(r,c_1,c_2,H)$ but to a subset
$S(r,c_1,c_2,A)_H$ of the $A$-Gieseker-semistable
sheaves consisting of those $E$ for which any subsheaf
$F \subset E$ having the same averaged Euler
characteristic \linebreak
with respect to $A$
$$p(F,A,n) = \frac{\chi(F \otimes A^n)}{rk(F)} =
\frac{\chi(E \otimes A^n)}{rk(E)} = p(E,A,n)$$
also satisfies
$$\frac{c_1(F)}{rk(F)} \cdot H \geq \frac{c_1(E)}{rk(E)}
\cdot H.$$
(Note that the direction of the inequality is opposite
from the one in the definition of
$H$-Gieseker-semistability.)  We denote the GIT quotient
$M(r,c_1,c_2,A)_H$.  This correspondence, while
seemingly technical and somewhat irregular, is the
important ingredient in the proof of the Main Theorem.

\S 3 begins the analysis of ``twists", which not only
gives us a second stratification of $\Delta - \{0\}$, built
on top of the previous one, characterizing the changes in
the set $S(r,c_1,c_2,H)$ as $H$ varies, but also tells us
how the (birational) map should be factorized into a
sequence of Thaddeus-type flips as in the Main Theorem.  A
torsion-free coherent sheaf $E$ is said to be ${\Cal
L}$-twisted $A$-Gieseker-semistable for ${\Cal L} \in
\text{Pic}(X) \otimes {\Bbb Q}$ if and only if for all $F
\subset E$
$$\frac{\chi(F \otimes {\Cal L} \otimes
A^n)}{rk(F)} \leq \frac{\chi(E \otimes {\Cal L} \otimes
A^n)}{rk(E)}$$
for $n >> 0$, where we compute the Euler characteristics
formally using the  \break
Riemann-Roch formula.  We denote the
set of such $E$ with given numerical data by \break
$S((r,c_1,c_2) \otimes {\Cal L},A)$, and in \S 4 and \S 5
we shall construct a projective moduli space
$M((r,c_1,c_2) \otimes {\Cal L},A)$ coarsely representing
the Seshadri equivalence classes of this set with respect
to the ${\Cal L}$-twisted $A$-Gieseker-stability.

Now consider the previous stratification $\Delta - \{0\}
= \coprod_s\Delta_s$, and denote by $V(\Delta_s)$ the
${\Bbb Q}$-span of $\Delta_s$ in $N^1(X)_{\Bbb Q}$.  Let
$\Delta_s$ and $\Delta_{s'}$ be two $d$-dimensional cells
such that $V(\Delta_s) = V(\Delta_{s'})$ and that they are
separated by a $d-1$-dimensional cell $W$.  Let $A$ be an
ample line bundle on $W$.  Then our analysis states that
there exists a stratification
$$V(\Delta_s) = \coprod_i L_i \coprod_j M_j$$
consisting of a finite number of hyperplanes $L_i$
parallel to $V(W)$ and connected components $M_j$ of
$V(\Delta_s) - \coprod_i L_i$ which determines the change
of the set $S((r,c_1,c_2) \otimes {\Cal L},A)$ as ${\Cal
L}$ varies in $V(\Delta_s)$.  Moreover, it follows that
$$\align
&S(r,c_1,c_2,A) = S(r,c_1,c_2,H) \text{\ \ for\ }A \in W
\text{\ and\ }H \in \Delta_s\\
\Longleftrightarrow&\\
&\overline{\Delta_s} \text{\ is\ contained\ in\ one\ of\ the\
strata\ }M_j\\
\endalign$$
and this gives a criterion, by induction on the dimension
$d$ of the cells, to determine the second stratification of
$\Delta - \{0\}$ characterizing the changes in the set
$S(r,c_1,c_2,H)$ as $H$ varies in $\Delta - \{0\}$.

Our aim then is to show that for a sequence of rational
ample line bundles ${\Cal L}_i \in L_i$ and ${\Cal M}_j
\in M_j$ there exists a sequence of Thaddeus-type flips:
$$
\CD
M((r,c_1,c_2) \otimes {\Cal M}_i,A) @.@.@.@.M((r,c_1,c_2)
\otimes {\Cal M}_{i+1},A)\\
@.\searrow @.@.\swarrow @.\\
@.@.M((r,c_1,c_2) \otimes {\Cal L}_i,A)@.@.\\
\endCD
$$
for $i = 0,1,\cdot\cdot\cdot,l$ arising from the Key GIT
Lemma.  This, together
 with the observation
$$\align
M(r,c_1,c_2,H) &= M((r,c_1,c_2) \otimes {\Cal M}_0,A)\\
M(r,c_1,c_2,H') &= M((r,c_1,c_2) \otimes {\Cal M}_{l+1},A)\\
M(r,c_1,c_2,A) &= M((r,c_1,c_2) \otimes {\Cal M}_i,A) \text{\
or\ }M((r,c_1,c_2) \otimes {\Cal L}_i,A) \text{\ for\ some\
}i,\\
\endalign$$
where $H \in \Delta_s$ and $H \in \Delta_{s'}$ (in fact,
we may take ${\Cal M}_0 = H^n$ and ${\Cal M}_{l+1} =
{H'}^{n'}$ for $n, n'$ sufficiently
 large) solves the problem of factorizing the
(birational) map between $M(r,c_1,c_2,H)$ and
$M(r,c_1,c_2,H')$, and thus proves the Main Theorem.

\S 4 is devoted to the Integral Case of the construction of
the flip.  When ${\Cal L}_i$ may be chosen to be an integral line
bundle, then $M((r,c_1,c_2) \otimes
{\Cal L}_i,A)$ is nothing but the classical moduli space
$M(r,c_{1_{{\Cal L}_i}},c_{2_{{\Cal L}_i}},A)$ where
$$\align
c_{1_{{\Cal L}_i}} &= c_1 + c_1({\Cal L}_i)\\
c_{2_{{\Cal L}_i}} &= c_2 + (r - 1)c_1 \cdot c_1({\Cal
L}_i) + \frac{r(r - 1)}{2}c_1({\Cal L}_i)^2.\\
\endalign$$
Moreover, we shall see that
$$\align
M((r,c_1,c_2) \otimes {\Cal M}_i,A) &\cong M(r,c_{1_{{\Cal
L}_i}},c_{2_{{\Cal L}_i}},A)_{H'}\\
M((r,c_1,c_2) \otimes {\Cal M}_{i+1},A) &\cong
M(r,c_{1_{{\Cal L}_i}},c_{2_{{\Cal L}_i}},A)_H.\\
\endalign$$
Therefore, the flip arises immediately from the construction in \S
2 and the Key GIT Lemma.

\S 5 deals with the more delicate case where ${\Cal L}_i$ is only a
RATIONAL line bundle.  Our strategy is to use Kawamata's
technique of finding some nice Galois cover $\phi:Y
\rightarrow X$ so that the set of ${\Cal L}_i$-twisted
$A$-Gieseker-semistable sheaves on $X$ correspond to the
 subset of ${\Cal L}_i^Y$-twisted $\phi^*A$-Gieseker-semistable sheaves
on $Y$, where ${\Cal L}_i^Y$ is now an INTEGRAL line bundle on
$Y$.  (Note that ${\Cal L}_i^Y$ is not exactly the pull
back $\phi^*{\Cal L}_i$ but ${\Cal L}_i^Y = \phi^*{\Cal
L}_i + \frac{1}{2}R$ where $R$ is the ramification
divisor of $\phi$.)  Then by looking at the appropriate locus
of the Quot scheme associated to $Y$ and taking the action of
the Galois group into consideration, we reduce the
construction of the flip to the Mumford-Thaddeus principle
and the Key GIT Lemma applied to the sheaves on $Y$ instead of $X$.

Finally starting from the $d$-cells of maximal dimension and
descending inductively on the dimension $d$, i.e., going
from $\Delta_s$ to $W$ and repeating, and moreover letting
$\Delta$ vary in $Amp(X)_{\Bbb Q}$, we can factorize all the
transformations among the moduli spaces $M(r,c_1,c_2,H)$ with
varying polarizations $H \in Amp(X)_{\Bbb Q}$ into sequences of
Thaddeus-type flips.  Our proof shows more generally that all the
transformations among the moduli spaces $M((r,c_1,c_2) \otimes {\Cal
L},H)$ with varying polarizations $H$ and rational twists ${\Cal L}$
can be factorized into sequences of Thaddeus-type flips.

\vskip.1in

We remark that instead of fixing the first chern class $c_1
\in N^1(X)_{\Bbb Q}$ we may fix the determinant of coherent
torsion-free sheaves to be $I \in \text{Pic}(X)$.  Then we obtain
$M(r,I,c_2,H)$
(resp. $M((r,I,c_2) \otimes {\Cal L},A)$) as a closed subscheme of
$M(r,c_1,c_2,H)$ (resp. $M((r,I,c_2) \otimes {\Cal L},A)$) where
$c_1(I) = c_1$.  The results above hold for these moduli
spaces with a fixed determinant without any change
(cf.[Gieseker-Li94]).

\vskip.1in

We would like to express our hearty thanks to V.
Alexeev, D. Abramovich, A. Bertram, A. Bruno, J. Harris, J. Koll\'ar,
S. Mori, M. Nakamaye, R. Pandharipande, D. Ruberman, M. Thaddeus and
K. Yoshioka for their invaluable suggestions and comments.  Their
help was more than crucial at several stages of the development
of the paper.  After this work was completed, R. Friedman and Z.
Qin kindly informed us that they had obtained results similar
to those in this paper in the case of rank 2 bundles on rational
surfaces with a more detailed analysis of flips
necessary for the computation of transition functions for
Donaldson polynomials [Friedman-Qin94].

\vskip.1in

\S 1. Some Basic Finiteness Results.

\vskip.1in

Let $X$ be a nonsingular projective surface over an
algebraically closed field $k$ of characteristic zero, $N^1(X)_{\Bbb
Q}$ the Neron-Severi group of $X$ tensored by ${\Bbb
Q}$, $Amp(X)_{\Bbb Q}$ the convex cone in
$N^1(X)_{\Bbb Q}$ generated by ample divisors.

\proclaim{Definition 1.1 (cf.[Gieseker77])} Let $H \in Amp(X)_{\Bbb Q}
\subset N^1(X)_{\Bbb Q}$ be a polarization.  A coherent
torsion-free sheaf $E$ of rank $rk(E)$ is
$H$-slope-semistable (resp. $H$-slope-stable) iff for all
subsheaves $F$ of $E$ (resp. $F
\underset\neq\to\subset E$)
$$\frac{c_1(F)}{rk(F)}\cdot H
\leq \frac{c_1(E)}{rk(E)} \cdot H\ (\text{\ resp.\ } < ).$$
$E$ is $H$-Gieseker-semistable (resp. $H$-Gieseker-stable)
iff for all coherent subsheaves $F$ of $E$
 (resp. $F \underset\neq\to\subset E$)
$$\align
p(F,H,n) := \frac{\chi(F \otimes H^n)}{rk(F)} &= \frac{1}{2}H^2n^2 +
(\frac{c_1(F)}{rk(F)}\cdot H - \frac{1}{2}K_X \cdot H)n\\
&+ \frac{1}{2}\frac{c_1(F)^2 - 2c_2(F) - c_1(F) \cdot
K_X}{rk(F)} + \chi({\Cal O}_X)\\
\leq p(E,H,n) = \frac{\chi(E \otimes H^n)}{rk(E)} &=
\frac{1}{2}H^2n^2 + (\frac{c_1(E)}{rk(E)}\cdot H -
\frac{1}{2}K_X \cdot H)n\\
&+ \frac{1}{2}\frac{c_1(E)^2 -
2c_2(E) - c_1(E) \cdot K_X}{rk(E)} + \chi({\Cal O}_X)\\
&(\text{resp.\ } p(F,H,n) < p(E,H,n))\\
&\text{\ for\ }n >> 0.
\endalign$$
We denote by $\mu(r,c_1,c_2,H)$ the set of all coherent torsion-free
sheaves of a fixed rank $r$, first and second chern
classes $c_1$ and $c_2$, which are $H$-slope-semistable, and use the
notation $S(r,c_1,c_2,H)$ for $H$-Gieseker-semistable ones.
\endproclaim

\vskip.1in

\proclaim{Definition 1.2} A set $T$ of coherent sheaves
on $X$ is bounded iff there is a scheme $S$ of finite
type over $k$ and a coherent sheaf ${\Cal F}$ over $S
\times X$ flat over $S$ so that for all $E \in T$ there
exists a closed point $s \in S$ such that the coherent
sheaf ${\Cal F}_s$ on the fiber of $S \times X$ over $s$
is isomorphic to $E$.
\endproclaim

\vskip.1in

\proclaim{Theorem 1.3} Let $\Delta$ be a convex cone
generated by a finite number of ample classes $H_1,
H_2, \cdot\cdot\cdot, H_l \in Amp(X)_{\Bbb Q} -
\{0\}$, i.e.,
$$\align
\Delta &= \{t_1H_1 + t_2H_2 + \cdot\cdot\cdot +
t_lH_l \in Amp(X)_{\Bbb Q};\\
&t_1, t_2, \cdot\cdot\cdot,t_l \in {\Bbb Q}, t_1 \geq
0, t_2 \geq 0, \cdot\cdot\cdot, t_l \geq 0\}\\
\endalign$$
Fix a numerical class $c_1 \in N^1(X)_{\Bbb Q}$, an
integer $c_2 \in {\Bbb Z}$ and another integer $r \in
{\Bbb N}$.

Then the set of all coherent torsion free sheaves $E$ of
$rk(E) = r$, $c_1(E) = c_1$ and $c_2(E) = c_2$ s.t. $E$ is
$H$-slope-semistable for some $H \in \Delta - \{0\}$, is
bounded.  That is to say, $\cup_{H \in \Delta -
\{0\}}\mu(r,c_1,c_2,H)$ is bounded.
\endproclaim

\demo{Proof of Theorem 1.3}\enddemo

We prove by induction on the rank $r$.

When $r = 1$, a coherent torsion-free sheaf $E$ of
$rk(E) = 1$, $c_1(E) = c_1$ and $c_2(E) = c_2$, is
isomorphic to ${\Cal L} \otimes {\Cal I}_Z$ where ${\Cal
I}_Z$ is the ideal sheaf defining an artinian scheme $Z$
of length $c_2$ and ${\Cal L}$ is a line bundle with
$c_1({\Cal L}) = c_1$.  Therefore, $S(1,c_1,c_2,H)$ is
parametrized by $\text{Pic}^{c_1}(X) \times
\text{Hilb}^{c_2}(X)$, where $\text{Pic}^{c_1}(X)$ is a
connected component of the Picard scheme of $X$ whose
corresponding line bundles have the same numerical
class $c_1$, and $\text{Hilb}^{c_2}(X)$ is the Hilbert
scheme parametrizing artinian schemes on $X$ of length
$c_2$.  Thus the family of these sheaves is bounded.
(Note that in the rank 1 case, the $H$-slope-semistability
is automatic.)

Next assuming the boundedness for the case of rank
$\leq r - 1$, we prove the boundedness for rank $= r$, by considering the
following two separate
subsets.

\vskip.1in

Case A: The subset consisting of those $E$ which are strictly
$H$-slope-semistable for some $H \in \Delta - \{0\}$,
i.e., there exists an exact sequence $$0 \rightarrow F
\rightarrow E \rightarrow G \rightarrow 0$$
where $F$ is a nonzero subsheaf of $E$ with $r(F) < rk(E) =
r$, $G$ is torsion free with the condition that
$$\frac{c_1(F)}{rk(F)} \cdot H = \frac{c_1(E)}{r} \cdot H =
\frac{c_1(G)}{rk(G)} \cdot H.$$
(If for the original choice of $F$ the cokernel $G$ is
not torsion-free, then we take the saturation of $F$ in
$E$ to make the cokernel torsion-free without changing the
first chern class of $F$ or $G$.)

\vskip.1in

For all such $E$ in this subset, we have
$$\align
c_1(F) &\equiv \frac{rk(F)}{r}c_1(E) + {\Cal L}\\
c_1(G) &\equiv \frac{rk(G)}{r}c_1(E) - {\Cal L}
\endalign$$
where ${\Cal L} \in \text{Pic}(X) \otimes \frac{1}{r}{\Bbb
Z}$ with ${\Cal L} \cdot H = 0$, and $\equiv$ denotes the
numerical equivalence.  Thus
$$\align
(*)\  c_2(E) &= c_1(F)c_1(G) + c_2(F) + c_2(G)\\
&= -{\Cal L}^2 + \frac{rk(G) - rk(F)}{r}c_1(E) \cdot {\Cal
L} + \frac{rk(F)rk(G)}{r^2}c_1(E)^2 \\
&\qquad + c_2(F) + c_2(G).
\endalign$$

We remark here that the Bogomolov inequality holds
not only for vector bundles but also for torsion-free
sheaves.

\proclaim{Lemma 1.4}  Let $E$ be a coherent torsion-free sheaf of rank $r$,
$H$-slope-semistable for some
$H \in Amp(X)_{\Bbb Q} - \{0\}$.  Then
$$c_2(E) \geq \frac{r - 1}{2r}c_1(E)^2.$$
\endproclaim

\demo{Proof of Lemma 1.4}\enddemo Let
$$0 \rightarrow E \rightarrow E^{**} \rightarrow Q
\rightarrow 0$$
be the inclusion of $E$ into its double dual $E^{**}$
with the cokernel $Q$ having a support on a finite
number of points.  Since $E^{**}$ is also
$H$-slope-semistable, by applying the Bogomolov inequality
to the vector bundle $E^{**}$ we conclude
$$c_2(E) =
c_2(E^{**}) + \text{length\ of\ } Q \geq c_2(E^{**}) \geq
\frac{r - 1}{2r}c_1(E^{**})^2 = \frac{r -
1}{2r}c_1(E)^2.$$

\vskip.1in

Since both $F$ and $G$ are $H$-slope-semistable, by Lemma 1.4 we have
$$\align
c_2(F) &\geq \frac{rk(F) - 1}{2rk(F)}c_1(F)^2\\
&= \frac{rk(F) - 1}{2rk(F)}\{{\Cal L}^2 +
2\frac{rk(F)}{r}c_1(E) \cdot {\Cal L} +
(\frac{rk(F)}{r})^2c_1(E)^2\},\\
c_2(G) &\geq \frac{rk(G) - 1}{2rk(G)}c_1(G)^2\\
&= \frac{rk(G) - 1}{2rk(G)}\{{\Cal L}^2 -
2\frac{rk(G)}{r}c_1(E) \cdot {\Cal L} +
(\frac{rk(G)}{r})^2c_1(E)^2\}.
\endalign$$

Plugging these into the formula for $c_2(E)$, we obtain
the inequality
$$\align
c_2 &= c_2(E)\\
&\geq (-1 + \frac{rk(F) - 1}{2rk(F)} +
\frac{rk(G) - 1}{2rk(G)}){\Cal L}^2\\
&+ 0 \cdot c_1(E) \cdot {\Cal
L} + \frac{(rk(F) + rk(G))(rk(F) + rk(G) -
1)}{2r^2}c_1(E)^2\\ &\geq (1 - h)(-{\Cal L}^2) + l,
\endalign$$
where $h$ and $l$ are constants depending only on $r$
and $c_1$
$$\align
1 > h &= \text{min}_{\{rk(F) + rk(G) = r\}}\{\frac{rk(F) - 1}{2rk(F)}
 + \frac{rk(G) -
1}{2rk(G)}\}\\ l &= \text{min}_{\{rk(F) + rk(G) = r\}}
\{\frac{(rk(F) + rk(G))(rk(F) + rk(G) -
1)}{2r^2}c_1(E)^2\}. \endalign$$
Therefore, we have
$$0 \leq -{\Cal L}^2 \leq \frac{c_2 - l}{1 - h}.$$

\proclaim{Lemma 1.5} Fix $N \in {\Bbb N}$.  The number of
lattice points $$\sharp\{x \in \text{Pic}(X) \otimes
\frac{1}{r}{\Bbb Z}; -x^2 \leq N, x \cdot H = 0 \text{\
for\ some\ }H \in \Delta\ - \{0\}\}$$
is finite.
\endproclaim

\demo{Proof of Lemma 1.5}\enddemo Observe that locally
with respect to $H \in \Delta - \{0\}$ we can find a
orthonormal basis of $N^1(X)_{\Bbb Q}$ consisting of
$e_1^H = \frac{H}{\sqrt{H \cdot H}}, e_2^H, \cdot\cdot\cdot,
e_{\text{dim}N^1(X)_{\Bbb Q}}^H$, which vary
continuously with respect to $H$ as seen, e.g, by the
method of Gram-Schmidt.  This implies that the set
$$\align
\{(H,x) \in \Delta \times N^1(X)_{\Bbb Q}&;-x^2 \leq
N, x \cdot H = 0,\\
&H = t_1H_1 + t_2H_2 + \cdot\cdot\cdot +
t_lH_l, t_1 + t_2 + \cdot\cdot\cdot +
t_l = 1\}
\endalign$$
is compact by the Hodge Index Theorem.
Thus its image in $N^1(X)_{\Bbb Q}$ is also compact, and
the number of the lattice points contained in it is finite.

\vskip.1in

By applying this lemma we conclude that the set of
the possible numerical classes of ${\Cal L}$ is finite,
which implies that both the set of the possible
numerical classes of $c_1(F)$ and that of the possible
numerical classes of $c_1(G)$ are finite.

Moreover, the finiteness of the numerical classes of
${\Cal L}$ implies that both $c_2(F)$ and $c_2(G)$ are
bounded from below, which follows from the Bogomolov
inequalities as above again.  This with the
equality ($*$) in turn implies that both $c_2(F)$ and
$c_2(G)$ have only a finite number of possibilities.

Now by inductional hypothesis on the rank and noting both $F$ and
$G$ are $H$-slope-semistable, we conclude that both $F$ and $G$ form
bounded families and so does $$E \in \text{Ext}^1(F, G).$$

We state the conclusion of this case in the following
form.

\proclaim{Proposition 1.6} The set of all torsion-free
coherent sheaves $F$ on $X$ such that there exists an
exact sequence
$$0 \rightarrow F \rightarrow E \rightarrow G
\rightarrow 0,$$
where $E \in \mu(r,c_1,c_2,H)$ for some $H \in \Delta - \{0\}$ ($r >
rk(F)$), and where the quotient $G$ is a coherent torsion-free sheaf,
satisfying $$\frac{c_1(F)}{rk(F)}\cdot H = \frac{c_1(E)}{r} \cdot
H,$$  is bounded, and so is the set of all such $G$.  In
particular, the sets of the pairs of the first and
second chern classes of such sheaves
$$\{(c_1(F),c_2(F))\}
\text{\ and\ } \{(c_1(G),c_2(G)\}$$
consist of a finite
number of elements. \endproclaim

\vskip.2in

We now deal with the second subset.

\vskip.1in

Case B: The complement of the previous subset, i.e., the
subset consisting of such $E$ that for
ANY $H\in \Delta - \{0\}$ is either $H$-slope-STABLE or
$H$-slope-UNSTABLE (equivalently, NOT
$H$-slope-semistable).

The condition of being $H$-slope-stable and that of being
$H$-slope-unstable are both open with respect to $H \in
\Delta - \{0\}$.  Since $E$ is $H$-slope-semistable for
some $H \in \Delta - \{0\}$, we conclude that $E$ must be
$H$-slope-stable for ALL $H \in \Delta - \{0\}$.  In
particular, $E$ is, say, $H_1$-slope-stable.  Thus $E$ is
$H_1$-Gieseker-semistable.  Therefore,
[Gieseker77,Corollary1.3] implies that the family of such
sheaves $E$
 is bounded.

This completes the proof of Theorem1.3.

\vskip.1in

\proclaim{Remark 1.7}\endproclaim

(1) The ``local-boundedness" as in Theorem 1.2 is the best we
can hope for IN GENERAL (though for some restricted
classes of surfaces the global-boundedness does hold as we
will see in (2)).  In fact, we can construct a counter-example
as follows to the statement that the set of all coherent
torsion free sheaves $E$ of fixed rank $rk(E) = r$, fixed
first and second chern classes $c_1(E) = c_1$
and $c_2(E) = c_2$, $H$-slope-semistable for
some $H \in Amp(X)_{\Bbb Q} - \{0\}$, is bounded.
(Actually our example presents an unbounded family
of coherent torsion free sheaves $E$ of fixed rank
$r$, fixed first and second chern
classes $c_1(E) = c_1$ and $c_2(E) = c_2$,
$H$-Gieseker-semistable for some $H \in
Amp(X)_{\Bbb Q} - \{0\}$.  Thus even if we restrict ourselves
to the Gieseker-semistable sheaves the global boundedness does not
hold.)

\vskip.1in

Take $X = E \times E$, where $E$ is a generic elliptic
curve (with no nontrivial automorphism fixing the origin other than
the involution) so that $\text{dim}_{\Bbb Q}N^1(X)_{\Bbb Q} = 3$.
Then $$\text{Pic}(X) = {\Bbb Z}E_1 \oplus {\Bbb Z}E_2 \oplus
{\Bbb Z}D (\text{\ in\ } N^1(X)_{\Bbb Q})$$
where $E_1 = E \times \{p\}, E_2 = \{q\} \times E$ and
$D$ is the diagonal.  Their intersection pairings are
$$\align
&E_1^2 = E_2^2 = D^2 = 0\\
&E_1 \cdot E_2 = E_1 \cdot D = E_2 \cdot D = 1\\
\endalign$$
and thus we have the intersection matrix
$$
\left(\matrix
0 &1 &1\\
1 &0 &1\\
1 &1 &0\\
\endmatrix\right).
$$
By diagonalizing the matrix, we have the corresponding
sublattice of rank 3 of $\text{Pic}(X)$ generated by
$$\matrix
U &= E_1 &+ E_2 &+ D\\
V &= E_1 &&- D\\
W &=\ \ \ &\ E_2 &- D.
\endmatrix$$
Their intersection pairings are
$$\align
U^2 &= 6, V^2 = W^2 = -2\\
U \cdot V &= V \cdot W = W \cdot U = 0.
\endalign$$
If we introduce the coordinate system $(u,v,w) = uU + vV +
wW$ in $N^1(X)_{\Bbb Q}$, then the ample cone
$Amp(X)_{\Bbb Q} - \{0\}$ is defined as the connected
component of $$\{(u,v,w) \in N^1(X)_{\Bbb Q};6u^2 - 2v^2 -
2w^2 > 0\}$$
containing $(1,0,0)$.

Note that if we take a line bundle ${\Cal L}$ whose
numerical class is given by
$${\Cal L} = uU + vV + wW \text{\ in\ }N^1(X)_{\Bbb Q},$$
then $c_2(E)$ of the rank 2 vector bundle
$$E = {\Cal L} \oplus {\Cal L}^{-1}$$
is given by
$$c_2(E) = - {\Cal L}^2 = - (6u^2 - 2v^2 - 2w^2).$$
Now we focus our attention to the plane $\{w = 0\}$.

\vskip2in

We recall an elementary lemma from number theory.

\proclaim{Lemma 1.8}
 Let $\omega \in {\Bbb R}_{> 0}$ be an irrational number.
Then there exist infinitely many rational numbers
$$\frac{q}{p}\text{\ where\ } p,q \in {\Bbb N}, g.c.d.(p,q)
= 1$$ such that
$$0 < \frac{q}{p} - \omega < \frac{1}{p^2}.$$
\endproclaim

We apply this lemma to $\omega = \sqrt{3}$.  Thus we have
infinitely many rational numbers
$$\frac{q}{p}\text{\ where\ } p,q \in {\Bbb N}, g.c.d.(p,q)
= 1$$ such that
$$0 < \frac{q}{p} - \sqrt{3} < \frac{1}{p^2}.$$
By excluding finitely many rational numbers we may assume
$$q < 2p.$$
Thus for a constant $c \geq 8$, we have
$$0 < q - \sqrt{3}p < \frac{1}{p} = \frac{8}{2(2p + 2p)}
< \frac{c}{2(q + \sqrt{3}p)},$$
which implies
$$0 < 2(q + \sqrt{3}p)(q - \sqrt{3}p) = - (6p^2 - 2q^2) <
c.$$
Therefore, there exists $c_2 \in {\Bbb N}$ with
$$0 < c_2 < c$$
such that there are infinitely many $(p,q) \in {\Bbb
N}^2$ with
$$- (6p^2 - 2q^2) = c_2.$$
Finally take line bundles ${\Cal L}_{(p,q)}$ s.t.
$${\Cal L}_{(p,q)} = pU + qV \text{\ in\ }N^1(X)_{\Bbb
Q}.$$
For the pairs $(p,q)$ with $q >> 0$, the point
$(p,q,0) \in N^1(X)_{\Bbb Q}$ is very close to the
boundary $\sqrt{3}u = v$ of the ample cone (intersected
with the plane $\{w = 0\}$).  Therefore, it follows
that there exists an (${\Bbb Q}$-)ample divisor
$H_{(p,q)}$ s.t.
$${\Cal L}_{(p,q)} \cdot H_{(p,q)} = 0.$$
(Actually $H_{(p,q)}$ gets very close inside of the
ample cone to the line $\sqrt{3}u = v$ as $q$
becomes large.)

The rank 2 vector bundles
$$E_{(p,q)} = {\Cal L}_{(p,q)} \oplus {\Cal
L}_{(p,q)}^{-1}$$
have the property
$$c_1(E_{(p,q)}) = 0, c_2(E_{(p,q)}) = c_2$$
and that $E_{(p,q)}$ is $H_{(p,q)}$-slope-semistable for
$H_{(p,q)} \in Amp(X)_{\Bbb Q} - \{0\}$.  Actually
$E_{(p,q)}$ is $H_{(p,q)}$-Gieseker-semistable since for
any $n$ we have
$$\align
\frac{1}{2}\chi(E_{(p,q)} \otimes H_{(p,q)}^n) &=
\chi({\Cal L}_{(p,q)} \otimes H_{(p,q)}^n) = \chi({\Cal
L}_{(p,q)}^{-1} \otimes H_{(p,q)}^n)\\
&= \frac{1}{2}n^2H_{(p,q)}^2 + \frac{1}{2}c_2.
\endalign$$

It is easy to see for any fixed ample line bundle $A$ on
$X$, the intersection number $A \cdot {\Cal L}_{(p,q)}$
is unbounded as $q$ becomes larger, thus the sheaf $E_{(p,q)}$ which
contains ${\Cal L}_{(p,q)}$ as a subsheaf is unbounded.

We remark also that what is essential in the example above
is the irrationality of the slope $\sqrt{3}$ and not the
circular shape of the cone.

\vskip.1in

(2) For some restricted classes of surfaces, the global
finiteness does hold.  We claim the global finiteness for
the class of surfaces whose nef cones are
rational and polyhedral (e.g. Del Pezzo surfaces).  To
see this, we prove a variant of Lemma 1.5.

\proclaim{Lemma 1.5'} Let $X$ be a nonsingular
projective surface whose nef cone is rational and
polyhedral, i.e., there exists a finite number of nef
line bundles $M_1, M_2, \cdot\cdot\cdot, M_m \in
\text{Pic}(X)$ such that $$\overline{Amp}(X)_{\Bbb R} =
{\Bbb R}_{\geq 0}M_1 + {\Bbb R}_{\geq 0}M_2 +
\cdot\cdot\cdot {\Bbb R}_{\geq 0}M_m.$$ Fix $N \in {\Bbb
N}$ and $r \in {\Bbb N}$.  Then the number of lattice
points
$$\sharp\{x \in \text{Pic}(X) \otimes \frac{1}{r}{\Bbb
Z};-x^2 \leq N, x \cdot H = 0 \text{\ for\ some\ } H \in
Amp(X)_{\Bbb Q} - \{0\}\}$$
is finite.
\endproclaim

\demo{Proof of Lemma 1.5'}\enddemo

To prove Lemma 1.5' we cannot use the argument of the proof of Lemma 1.5
directly.  But essentially the only thing that may give
rise to the infinite number of lattice points is the
existence of some ``irrational" edge on the boundary of
the ample cone as in the previous counterexample, which is
excluded by the assumption.  Since after the application of Hodge index Theorem
 the proof is just an analysis of lattice points bounded by some quadratic
hypersurface
 (and since we will not use Lemma 1.5' or Theorem 1.3' in the arguments of
later
chapters), we leave the details of the proof to the reader.

Once we have the finiteness of the lattice points above,
the rest of the argument goes through without change.
Thus we have

\proclaim{Theorem 1.3'} Let $X$ be a nonsingular
projective surface whose nef cone is rational and
polyhedral (e.g., a Del Pezzo surface, a relatively
minimal ruled surface over a smooth curve (cf.[CKM88]) or a special kind of K3
surface (cf. [Kov\'acs94])).  Fix
a numerical class $c_1 \in N^1(X)_{\Bbb Q}$, an integer
$c_2 \in {\Bbb Z}$ and another integer $r \in {\Bbb
N}$.  Then the set of all coherent torsion-free sheaves
$E$ of $rk(E) = r$, $c_1(E) = c_1$ and $c_2(E) = c_2$ s.t.
$E$ is $H$-slope-semistable for some $H \in Amp(X)_{\Bbb
Q} - \{0\}$, is bounded.
\endproclaim

We also have the boundedness for Gieseker-semistable sheaves.

\proclaim{Corollary 1.9} Under the same assumption as
in Theorem 1.3 (resp. Theorem 1.3'), the set of all
coherent torsion free sheaves $E$ of $rk(E) = r$, $c_1(E) =
c_1$ and $c_2(E) = c_2$, $H$-Gieseker-semistable for some
$H \in \Delta - \{0\}$ (resp. $H \in$\linebreak
$ Amp(X)_{\Bbb Q} -
\{0\}$), is bounded.  That is to say, $\cup_{H \in \Delta
- \{0\}}S(r,c_1,c_2,H)$ (resp. \linebreak
$\cup_{H \in Amp(X)_{\Bbb
Q} - \{0\}}S(r,c_1,c_2,H)$) is bounded. \endproclaim

\demo{Proof of Corollary 1.9}\enddemo

If $E$ is $H$-Gieseker-semistable, then $E$ is
$H$-slope-semistable.  That is to say, the set of
$H$-Gieseker-semistable sheaves is a subset of
$H$-slope-semistable sheaves.  Now the assertion of the
corollary is immediate from that of Theorem 1.3 (resp.
Theorem 1.3').

We will use the boundedness results of this section in \S 3 to
determine the stratification of $\Delta - \{0\}$ which describes the
change of the set $\mu(r,c_1,c_2,H)$ as $H$ varies in $\Delta -
\{0\}$ and that of the set $S(r,c_1,c_2,H)$.

\vskip.1in

\S 2. Key GIT Lemma after Simpson

\vskip.1in

In \S 2, we focus our attention on the key GIT lemma,
modifying the ideas of [Simpson92], which gives us the main
machinery needed to reduce the problem of factorizing the
(birational) transformations among various moduli spaces to
the Mumford-Thaddeus principle.  We use the same notation
as in \S 1.  We shall present the entire argument;
while the reader may find this a bit repetitive of
[Simpson92], we believe the
necessary modifications  too subtle and too many
to just leave to reference.

\vskip.1in

Step 1:Embedding into a Grothendieck Quot scheme.

\vskip.1in

Let $\cup_{H \in \Delta - \{0\}}S(r,c_1,c_2,H)$ be the
set of all coherent torsion-free sheaves of rank
$r$, fixed first and second chern classes $c_1$
and $c_2$, $H$-Gieseker-semistable for SOME $H \in \Delta
- \{0\}$.  Then the result of \S 1 tells us that the set
$\cup_{H \in \Delta - \{0\}}S(r,c_1,c_2,H)$ is bounded,
i.e., there exists a scheme $S$ of finite type over $k$ and
coherent sheaf ${\Cal F}$ on $S \times X$ flat over $S$
such that for any $E \in \cup_{H \in \Delta -
\{0\}}S(r,c_1,c_2,H)$ there is $s \in S$ with $E \cong
{\Cal F}_s$.  Let $\pi_1:S \times X \rightarrow S$ and
$\pi_2:S \times X \rightarrow X$ be the first and second
projections respectively.

\vskip.1in

Let $A$ be an ample divisor on $X$ with $A \in \Delta -
\{0\}$.

\vskip.1in

Then there exists $a_0 \in {\Bbb N}$ such that for all $a
\geq a_0$
$${\pi_1}^*{\pi_1}_*({\Cal F} \otimes {\pi_2
}^*A^a) \rightarrow {\Cal F} \otimes {\pi_2}^*A^a$$
is surjective, and
$$R^i{\pi_1}_*({\Cal F} \otimes {\pi_2
}^*A^a) = 0 \text{\ for\ } i > 0.$$

Furthermore, for all $s \in S$ we have that
$$R^i{\pi_1}_*({\Cal F} \otimes {\pi_2}^*A^a) \otimes
k(s) \rightarrow H^i(X,{\Cal F}_s \otimes A^a)$$
is an isomorphism for $i \geq 0$ and that ${\pi_1
}_*({\Cal F} \otimes {\pi_2}^*A^a)$ is locally free of
rank $l = \chi_{r,c_1,c_2}(A^a)$ (cf. [Hartshorne,
Cohomology and Base Change, Theorem12.11]), where we
introduce the notation
$$\align
\chi_{r,c_1,c_2}(A^a) &= \frac{r}{2}A^2a^2 +
(c_1 \cdot A - \frac{r}{2}K_X \cdot A)a\\
&+ \frac{1}{2}(c_1^2 - 2c_2 - c_1 \cdot K_X) +
r\chi({\Cal O}_X).
\endalign$$
We may assume $S$ is a finite disjoint union of affine
schemes $S_{\alpha}$ such that
$${\pi_1}_*({\Cal F} \otimes {\pi_2}^*A^a)|_{S_{\alpha}}
\cong {\Cal O}_{S_{\alpha}}^{\oplus l}$$
and thus we have a surjection
$${\Cal O}_{S \times X}^{\oplus l} \cong
{\pi_1}^*{\pi_1}_*({\Cal F} \otimes {\pi_2}^*A^a)
\rightarrow {\Cal F} \otimes {\pi_2}^*A^a.$$
Therefore, we obtain a morphism
$$\phi:S \rightarrow Quot({\Cal O}_X^{\oplus
l}/\chi_{r,c_1',c_2'})$$
such that for all $s \in S$ we have
$$({\Cal F} \otimes {\pi_2}^*A^a)_s \cong
(Univ)_{\phi(s)},$$
where $Quot({\Cal O}_X^{\oplus
l}/\chi_{r,c_1',c_2'})$ (Later we will use the
abbreviation $Quot$.) is the Grothendieck's Quot scheme
parametrizing the quotients of ${\Cal O}_X^{\oplus l}$
whose Hilbert polynomial with respect to an ample line
bundle $H$ is given by $\chi_{r,c_1',c_2'}(H^m)$, and
$$\align c_1' &= c_1 + rac_1(A)\\
c_2' &= c_2 + (r - 1)ac_1 \cdot c_1(A) + \frac{r(r -
1)}{2}a^2c_1(A)^2.
\endalign$$

\vskip.1in

Step 2:Inducing the linearized polarizations onto the
Quot scheme.

\vskip.1in

We take an ample line bundle $H$ on $X$.

\vskip.1in

Let $\pi_1:Quot
\times X \rightarrow Quot$ and $\pi_2:Quot \times X
\rightarrow X$ be the first and second projections.  We
have the universal quotient sheaf $Univ$ over $Quot
\times X$
$$0 \rightarrow Ker \rightarrow {\Cal O}_{Quot
\times X}^{\oplus l} \rightarrow Univ \rightarrow 0.$$

There exists $M_H \in {\Bbb N}$ such that for all $m
\geq M_H$ we have
$$\align
R^i{\pi_1}_*(Univ \otimes {\pi_2}^*H^m) &= 0\\
R^i{\pi_1}_*({\Cal O}_{Quot \times X}^{\oplus l} \otimes
{\pi_2}^*H^m) &= 0\\
R^i{\pi_1}_*(Ker \otimes {\pi_2}^*H^m) &= 0\\
&\text{\ for\ }i > 0,
\endalign$$
and for all closed points $q \in Quot$
$$\align
R^i{\pi_1}_*(Univ \otimes {\pi_2}^*H^m) \otimes k(q)
&\rightarrow H^i(X,(Univ)_q \otimes H^m)\\
R^i{\pi_1}_*({\Cal O}_{Quot \times X}^{\oplus l} \otimes
{\pi_2}^*H^m) \otimes k(q) &\rightarrow H^i(X,{\Cal
O}_X^{\oplus l} \otimes H^m)\\
R^i{\pi_1}_*(Ker \otimes {\pi_2}^*H^m) \otimes k(q)
&\rightarrow H^i(X,(Ker)_q \otimes H^m)\\
\endalign$$
are all isomorphisms for $i \geq 0$.

Furthermore,
$${\pi_1}_*(Univ \otimes {\pi_2}^*H^m)$$
is locally free of rank $R_m = \chi_{r,c_1',c_2'}(H^m)$,
$${\pi_1}_*({\Cal O}_{Quot \times X}^{\oplus l} \otimes
{\pi_2}^*H^m) \cong {\Cal O}_{Quot}^{\oplus l} \otimes_k
H^0(X,H^m)$$
is locally free and so is
$${\pi_1}_*(Ker \otimes {\pi_2}^*H^m).$$
Therefore we have the surjection
$${\Cal O}_{Quot}^{\oplus l} \otimes_k H^0(X,H^m) \cong
{\pi_1}_*({\Cal O}_{Quot \times X}^{\oplus l} \otimes
{\pi_2}^*H^m) \rightarrow {\pi_1}_*(Univ \otimes
{\pi_2}^*H^m) \rightarrow 0$$
with ${\pi_1}_*(Univ \otimes
{\pi_2}^*H^m)$ being locally free, which induces a
morphism from the Quot scheme into the Grassmannian of the
$R_m$-dimensional quotients of $k^{\oplus l} \otimes_k
H^0(X,H^m)$  $$\phi_m:Quot \rightarrow Grass(k^{\oplus l}
\otimes_k H^0(X,H^m),R_m).$$
By taking $M_H$ sufficiently large, we may assume that
$\phi_m$ is an embedding.

The algebraic group $SL(l)$ naturally acts both on $Quot$
and $Grass(k^{\oplus l} \otimes_k H^0(X,H^m),R_m)$, and
$\phi_m$ is $SL(l)$-linear ($SL(l)$-equivariant).  There
is a natural $SL(l)$-linear embedding of
$Grass(k^{\oplus l} \otimes_k H^0(X,H^m),R_m)$ by
 the Pl\"ucker embedding, which gives an
$SL(l)$-linearization on the ample line bundle on $Quot$
$$det\{{\pi_1}_*(Univ \otimes {\pi_2}^*H^m)\}.$$ Thus Mumford's
notion of semistable (stable) points w.r.t. this
linearized ample bundle is defined on $Quot$.

The Hilbert-Mumford numerical criterion for
semistability (stability) for the Grassmannian $Grass(V
\otimes W,R)$ of the $R$-dimensional quotients of the
vector space $V \otimes W$ under the action of $SL(V)$ with
the linearization induced by the Pl\"ucker embedding is
given by the following lemma.

\proclaim{Lemma 2.1} A point $p:V \otimes W
\rightarrow U \rightarrow 0$ in $Grass(V \otimes W,R)$
is semistable (resp. stable) for the action of $SL(V)$
and the linearization induced by the Pl\"ucker
embedding if and only if for all nonzero proper
subspaces $0 \neq L \underset\neq\to\subset V$ we have
$p(L \otimes W) \neq 0$ and
$$\frac{\text{dim}L}{\text{dim}\ p(L \otimes W)} \leq
\frac{\text{dim}V}{\text{dim}U}\ (\text{resp.\ }<).$$
\endproclaim

See [Simpson92,Proposition 1.14] and
for the proof [MF82,Proposition 4.3].

\vskip.1in

\noindent
[Simpson92,Lemma1.15] also gives a criterion for a point $q \in
Quot$ to be semistable.

\proclaim{Lemma 2.2} Given the situation as in Step
1 and Step 2 above.  Then there exists $M_H \in {\Bbb N}$
such that for all $m \geq M_H$ the following holds: Suppose
$q:{\Cal O}_X^{\oplus l} \rightarrow E \otimes A^a
\rightarrow 0$ is a point in $Quot$ such that
$$k^{\oplus l} \cong H^0({\Cal O}_X^{\oplus l})
\rightarrow H^0(E \otimes A^a)$$
is an isomorphism, and that
$$\frac{h^0(F_L \otimes A^a)}{\chi(F_L \otimes A^a \otimes
H^m)} \leq \frac{h^0(E \otimes A^a)}{\chi(E \otimes A^a
\otimes H^m)} = \frac{l}{\chi_{r,c_1',c_2'}(H^m)}
(\text{resp.\ }<)$$
for all subsheaves $F_L
\subset E$ where $F_L \otimes A^a$ is generated by a
nonzero 9resp. nonzero prper) linear subspace $L$ of $H^0(E
\otimes A^a)$.  Then $q$ is semistable (resp. stable) w.r.t. the
action of $SL(V)$ and the linearization induced from the Pl\"ucker
embedding of $Grass(k^{\oplus l} \otimes_k H^0(X,H^m),R_m)$.
\endproclaim

\demo{Proof of Lemma 2.2}\enddemo

For all nonzero linear subspaces $0 \neq L \subset V =
k^{\oplus l} = H^0({\Cal O}_X^{\oplus l}) \cong H^0(E
\otimes A^a)$, let $F_L$ be the subsheaf of $E$ such that
$F_L \otimes A^a$ is generated by $L$.  Note
$$\text{dim}L \leq h^0(F_L \otimes A^a).$$
$Quot$ is embedded into $Grass(k^{\oplus l} \otimes_k
H^0(X,H^m),R_m) = Grass(V \otimes W,R)$ (if we take
$M_H$ large enough), which in turn is embedded into a
projective space by the Pl\"ucker embedding.

Remark that once $a$ is fixed the family of subsheaves
$F_L \otimes A^a$ of $E \otimes A^a = (Univ)_q$ generated
by a linear subspace $L$ as $q$ varies among all points
in $Quot$ and $L$ varies among all subspaces of
$H^0({\Cal O}_X^{\oplus l})$ is bounded; similarly,
the family of such kernels $K_L$ of the exact
sequences   $$0 \rightarrow K_L \rightarrow L
\otimes_k {\Cal O}_X \rightarrow F_L \otimes A^a \rightarrow
0$$ is bounded.

Now take $M_H \in {\Bbb N}$ so that for all $m \geq M_H$
we have
$$h^0(F_L \otimes A^a) \leq h^0(F_L \otimes A^a
\otimes H^m) = \chi(F_L \otimes A^a \otimes H^m)$$
and
$$h^1(K_L \otimes H^m) = 0$$
for all such $F_L$ and $K_L$.
Thus from the exact sequence
$$0 \rightarrow K_L \otimes H^m \rightarrow L \otimes_k H^m
\rightarrow F_L \otimes A^a \otimes H^m \rightarrow 0$$
we have
$$L \otimes_k H^0(X,H^m) \rightarrow H^0(F_L \otimes A^a
\otimes H^m) \rightarrow H^1(K_L \otimes H^m) = 0,$$
which implies
$$\text{dim}\ p(L \otimes W) = h^0(F_L \otimes A^a \otimes
H^m) = \chi(F_L \otimes A^a \otimes H^m).$$
Therefore, we conclude
$$\align
\frac{\text{dim}L}{\text{dim}\ p(L \otimes W)} &\leq
\frac{h^0(F_L \otimes A^a)}{\text{dim}\ p(L \otimes W)}\\
&= \frac{h^0(F_L \otimes A^a)}{\chi(F_L \otimes A^a
\otimes H^m)}\\
&\leq \frac{h^0(E \otimes A^a)}{\chi(E \otimes A^a
\otimes H^m)}\\
&= \frac{\text{dim}V}{\text{dim}U}.\\
\endalign$$
(resp. If $L$ is a proper subspace, then either $F_L = E$
and the first inequality is strict $<$, or $F_L$ is a
proper subsheaf of $E$ and the second inequality is
strict $<$ by assumption.)

Thus by Lemma 2.1 $q$ is semistable (resp. stable).

\vskip.1in

The following lemma, a variant of [Simpson92,Lemma1.16],
which characterizes a semistable point $q:{\Cal
O}_X^{\oplus l} \rightarrow E \otimes A^a \rightarrow 0$
plays a crucial role with the previous lemma in establishing
the key GIT lemma of this section.

\proclaim{Lemma 2.3} For a fixed $a \in {\Bbb N} (\geq
a_0)$ and an ample line bundle $H$ on $X$, there exists
$M_H \in {\Bbb N}$ s.t. for all $m \geq M_H$, the
following holds: If a point
$$q:{\Cal O}_X^{\oplus l} \rightarrow E \otimes A^a
\rightarrow 0$$
is semistable with respect to the
action of $SL(l)$ and the linearization induced from
the Pl\"ucker embedding of $Grass(k^{\oplus l} \otimes_k
H^0(X,H^m),R_m)$, then the natural homomorphism
$$k^{\oplus l} = H^0({\Cal O}_X^{\oplus l}) \rightarrow
H^0(E \otimes A^a)$$
is injective, and for any nonzero quotient
$$E \otimes A^a
\rightarrow G \otimes A^a \rightarrow 0$$
with $rk(G) > 0$, we have
$$\align
\frac{h^0(G \otimes A^a)}{rk(G)} &\geq \frac{\chi(E
\otimes A^a)}{rk(E)} = \frac{l}{r}.
\\
\endalign$$
Moreover, suppose that ${\Cal G}$ over $T \times X$ is a
bounded family of coherent sheaves on $X$ (independent of
$a$).  Then there exists $a_2 \in {\Bbb N} (\geq a_0)$ such
that for all $a \geq a_2$ and an ample line bundle $H$ on
$X$, there exists $M_H \in {\Bbb N}$ such that for all
$m \geq M_H$ the following holds: If a point   $$q:{\Cal
O}_X^{\oplus l} \rightarrow E \otimes A^a \rightarrow 0$$
is semistable with respect to the
action of $SL(l)$ and the linearization induced from
the Pl\"ucker embedding of $Grass(k^{\oplus l} \otimes_k
H^0(X,H^m),R_m)$ and if the natural homomorphism
$$k^{\oplus l} = H^0({\Cal O}_X^{\oplus l}) \rightarrow
H^0(E \otimes A^a)$$
is an isomorphism, then for any nonzero quotient
$$E \otimes A^a \rightarrow G \otimes A^a \rightarrow 0$$
with $G \cong {\Cal G}_t$ for some $t \in T$, we have
$$
\align\frac{h^0(G \otimes A^a)}{\chi(G \otimes A^a \otimes
H^m)} &\geq \frac{\chi(E \otimes A^a)}{\chi(E \otimes A^a
\otimes H^m)}.\\
\endalign$$
\endproclaim

\demo{Proof of Lemma 2.3}\enddemo

Suppose a point $q:{\Cal O}_X^{\oplus l} \rightarrow E
\otimes A^a \rightarrow 0$ is semistable.  Then Lemma
2.1 tells us that for any nonzero linear subspace $0
\neq L \subset V = H^0({\Cal O}_X^{\oplus l})$, $p$
being the map given by
$$\align
p:L \otimes W (= &H^0(H^m)) \rightarrow U (= H^0(E \otimes
A^a \otimes H^m))\\
&\searrow\ \ \ \ \ \ \ \ \ \ \ \ \ \ \ \ \ \ \ \
\ \nearrow\\
&H^0(E \otimes A^a) \otimes H^0(H^m)\\
\endalign$$
we have $p(L \otimes W) \neq 0$.  Thus the homomorphism
$L \rightarrow H^0(E \otimes A^a)$ is nonzero.  Since
$L$ is an arbitrary nonzero linear subspace, we have
the desired injectivity of $$k^{\oplus l} = H^0({\Cal
O}^{\oplus l}) \rightarrow H^0(E \otimes A^a).$$

Suppose
$$E \otimes A^a \rightarrow G \otimes A^a \rightarrow 0$$
is a nonzero quotient with $rk(G) > 0$ and
$$
\align\frac{h^0(G \otimes A^a)}{rk(G)} &< \frac{\chi(E
\otimes A^a)}{rk(E)} = \frac{l}{r}.\\
\endalign$$
Let $L (\subset V)$ be the kernel of the homomorphism
$$V = H^0({\Cal O}_X^{\oplus l}) (\rightarrow H^0(E
\otimes A^a)) \rightarrow H^0(G \otimes A^a),$$
which is nonzero since we have the strict inequality
above and $rk(G) \leq rk(E)$ and hence
$$h^0(G \otimes A^a) < l = h^0({\Cal O}_X^{\oplus l}).$$
Let $F_L$ be the subsheaf of $E$ such that $F_L \otimes
A^a$ is generated by (the image of) $L$ in $H^0(E \otimes
A^a)$.

Note that from the exact sequence
$$0 \rightarrow K_G \otimes A^a \rightarrow E \otimes
A^a \rightarrow G \otimes A^a \rightarrow 0$$
and the inclusion
$$F_L \otimes A^a \hookrightarrow K_G \otimes A^a,$$
it follows
$$(*)\ rk(F_L) + rk(G) \leq rk(E).$$
Moreover, we have
$$(**)\ \text{dim}L \geq l - h^0(G \otimes A^a).$$
Therefore, the strict inequality of the
assumption together with $(*)$ and $(**)$ implies
$$\align
\frac{\text{dim}L}{rk(F_L)} &> \frac{l}{rk(E)}.\\
\endalign$$
Thus we obtain
$$
\align(\heartsuit)\ \frac{\text{dim}L}{\chi(F_L \otimes A^a
\otimes H^m)} &> \frac{l}{\chi(E \otimes A^a \otimes
H^m)}\\
\endalign$$
for $m >> 0$.

Note that once $a$ is fixed, such $F_L$ (as well as
$F_L \otimes A^a$) ranges over a bounded family and
thus we have a finite number of possibilities for
$\chi(F_L \otimes A^a \otimes H^m)$ as polynomials in
$m$.

Take $M_H \in {\Bbb N}$ such that for all $m \geq M_H$
we have
$$\align
h^0(F_L \otimes A^a \otimes H^m) &= \chi(F_L \otimes
A^a \otimes H^m)\\
h^0(E \otimes A^a \otimes H^m) &= \chi(E \otimes
A^a \otimes H^m),
\endalign$$
$$h^1(K_L \otimes H^m) = 0$$
where
$$0 \rightarrow K_L \rightarrow L \otimes_k {\Cal O}_X
\rightarrow F_L \otimes A^a \rightarrow 0,$$
and that $(\heartsuit)$ holds if we have such $G$.

Then
 we have$$\align
\text{dim}\ p(L \otimes W) &= h^0(F_L \otimes A^a \otimes
H^m)\\
&= \chi(F_L \otimes A^a \otimes H^m),
\endalign$$
which implies
$$\align
\frac{\text{dim}L}{\text{dim}\ p(L \otimes W)} &>
\frac{l}{h^0(E \otimes A^a \otimes H^m)}.\\
\endalign$$
This contradicts the criterion for semistability of $q$
given in Lemma 2.1.

\vskip.1in

Now we proceed to the proof of the second
 half of the lemma.

Note first that since the family of such $G$ is bounded
by assumption, so is the family of kernels $K_G$ of
the sequence
$$0 \rightarrow K_G \rightarrow E \rightarrow G \rightarrow
0.$$
Therefore, there exists $a_2 \in {\Bbb N} (\geq a_0)$
such that for all $a \geq a_2$ and exact sequences
$$0 \rightarrow K_G \otimes A^a \rightarrow E \otimes A^a
\rightarrow G \otimes A^a \rightarrow 0,$$
we have the exactness of
$$0 \rightarrow H^0(K_G \otimes A^a) \rightarrow H^0(E
\otimes A^a) \rightarrow H^0(G \otimes A^a) \rightarrow
0$$
by the virtue of
$$H^1(K_G \otimes A^a) = 0,$$
and that $K_G \otimes A^a$ is generated by the global
sections $H^0(K_G \otimes A^a)$.

Since we assume
$$k^{\oplus l} \cong H^0({\Cal O}_X^{\oplus l})
\rightarrow H^0(E \otimes A^a)$$
is an isomorphism, we have
$$L = H^0(K_G \otimes A^a) \text{\ and\ }F_L = K_G$$
in the notation of the first half of our proof.

Therefore, we have
$$(*')\ \chi(F_L \otimes A^a \otimes H^m) + \chi(G \otimes
A^a \otimes H^m) = \chi(E \otimes A^a \otimes H^m)$$
and
$$(**')\ \text{dim}L = l - h^0(G \otimes A^a).$$
Take $M_H \in {\Bbb N}$ so that for all $m \geq M_H$ we
have
$$\align
\chi(F_L \otimes A^a \otimes H^m) &= h^0(F_L \otimes A^a
\otimes H^m) > 0\\
\chi(G \otimes A^a \otimes H^m) &= h^0(G \otimes A^a
\otimes H^m) > 0\\
\chi(E \otimes A^a \otimes H^m)& = h^0(E \otimes A^a
\otimes H^m) > 0,\\
\endalign$$
and that
$$h^1(K_L \otimes H^m) = 0$$
where
$$0 \rightarrow K_L \rightarrow L \otimes_k {\Cal O}_X
\rightarrow F_L \otimes A^a \rightarrow 0.$$
Suppose
$$\align
\frac{h^0(G \otimes A^a)}{\chi(G \otimes A^a \otimes
H^m)} &< \frac{\chi(E \otimes A^a)}{\chi(E \otimes A^a
\otimes H^m)}
 = \frac{l}{\chi(E \otimes A^a \otimes H^m)}.
\\
\endalign$$
This together with $(*')$ and $(**')$ implies
$$\align
(\heartsuit)\ \frac{\text{dim}L}{\chi(F_L \otimes A^a
\otimes H^m)} &> \frac{\chi(E \otimes A^a)}{\chi(E \otimes
A^a \otimes H^m)}.\\
\endalign$$
Now
$$h^1(K_L \otimes H^m) = 0$$
implies
$$\align
\text{dim}\ p(L \otimes W) &= h^0(F_L \otimes A^a \otimes
H^m)\\
&= \chi(F_L \otimes A^a \otimes H^m).
\endalign$$
Thus $(\heartsuit)$ leads to
$$\align
\frac{\text{dim}L}{\text{dim}\ p(L \otimes W)} &>
 \frac{l}{h^0(E \otimes A^a \otimes H^m)},\\
\endalign$$
which again contradicts the criterion of semistability of
$q$ given in Lemma 2.1.

Ideally and quite naively one would
expect after having embedded all the sheaves in
$$\cup_{H
\in \Delta - \{0\}}S(r,c_1,c_2,H)$$
into the Quot scheme
$Quot$ by taking a sufficiently high multiple of an ample line bundle $A \in
\Delta - \{0\}$, the semistable
points with respect to the $SL(l)$-action and the
linearization induced from the Pl\"ucker embedding
twisting by a high multiple of an ample line bundle $H \in Amp(X)_{\Bbb Q} -
\{0\}$ would correspond to the
$H$-Gieseker-semistable sheaves.  In fact in the classical
setting with the fixed polarization $A = H$ this
is indeed the case.  But in our setting with changing
polarizations, the situation becomes more subtle.  Steps 1
and 2 are not independent but are actually quite intertwined,
and the naive expectation as above turns out to be false.
But the study of the change of the sets of the semistable
points as we change the second polarization $H$ gives the
vital information to analyze the transformation among the
various moduli spaces.

\proclaim{Key GIT Lemma 2.4}  Let $A$ be a VERY ample
line bundle with $A \in \Delta - \{0\}$.  Then there exists
$a_A \in {\Bbb N}$ such that for all $a \geq a_A$ the
following holds:

\ \ (i) All the sheaves in $\cup_{H
\in \Delta - \{0\}}S(r,c_1,c_2,H)$ can be embedded into
a Quot scheme by taking a multiple $A^a$ of $A$,
i.e., for all $E \in \cup_{H
\in \Delta - \{0\}}S(r,c_1,c_2,H)$ there exists a point
$$q \in Quot({\Cal O}_X^{\oplus l}/\chi_{r,c_1',c_2'})$$
s.t.
$$E \otimes A^a \cong (Univ)_q$$
where
$$\align
l &= \chi_{r,c_1,c_2}(A^a)\\
c_1' &= c_1 + rac_1(A)\\
c_2' &= c_2 + (r - 1)c_1 \cdot ac_1(A) + \frac{r(r -
1)}{2} \cdot a^2c_1(A)^2.
\endalign$$

\ \ (ii) We denote by $Q$ the closure in $Quot({\Cal O}_X^{\oplus
l}/\chi_{r,c_1',c_2'})$ of the set of points $q$ such
that $(Univ)_q$ is torsion free with $c_1((Univ)_q) = c_1'$ and
$c_2((Univ)_q) = c_2'$.  Note that $Q$ is $SL(l)$-invariant.

For any ample line bundle $H \in Amp(X)_{\Bbb Q} -
\{0\}$, there exists $M_H \in {\Bbb N}$ such that for all
$m \geq M_H$ the following equivalence holds:

A point $q \in Q$ is semistable with respect to the action
of $SL(l)$ and the linearization induced from the Pl\"ucker
embedding of $Grass(k^{\oplus l} \otimes H^0(X,H^m),R_m)$

\vskip.1in

if and only if

\vskip.1in

$(Univ)_q = E \otimes A^a$ is a coherent
torsion free sheaf of rank $rk(E) = r$, $c_1(E) = c_1$,
$c_2(E) = c_2$ such that $E$ is $A$-Gieseker-semistable
and that for all subsheaves $F \subset E$
having the same averaged Euler characteristics
$p(F,A,a) = p(E,A,a)$ (as polynomials in $a$) we have
$$\frac{c_1(F)}{rk(F)} \cdot H \geq
\frac{c_1(E)}{rk(E)} \cdot H,$$
i.e.,
$$E \in S(r,c_1,c_2,A)_H,$$
and furthermore the natural homomorphism
$$k^{\oplus l} = H^0({\Cal O}_X^{\oplus l}) \rightarrow
H^0((Univ)_q)$$
is an isomorphism.
\endproclaim

\demo{Proof of Lemma 2.4}\enddemo

Remark first that the assumption $A$ being VERY ample is
only for the sake of simplicity.  For
an arbitrary ample line bundle $A$, one just has to take
a multiple to make it very ample and the corresponding
statements hold.

\ \ (i) This is the direct consequence
 of the finiteness results of \S 1, and the construction
is explained at the beginning of \S 2.

\ \ (ii) We start by proving the following preliminary
claim (cf. [Simpson92, Lemma1.18]).

\proclaim{Claim 2.5} There exists $a_1 \in {\Bbb
N}$ such that for all $a \geq a_1$ the following
holds:

For all subsheaves $F$ of any $E \in S(r,c_1,c_2,A)$ we
have $$\frac{h^0(F \otimes A^a)}{rk(F)} \leq \frac{h^0(E
\otimes A^a)}{rk(E)} = \frac{\chi(E \otimes
A^a)}{rk(E)}$$
and the equality holds if and only if
$$\frac{\chi(F \otimes A^a \otimes A^m)}{rk(F)} =
\frac{\chi(E \otimes A^a \otimes A^m)}{rk(F)}$$
for all $m$, i.e.,
$$p(F,A,a) = p(E,A,a)$$
as polynomials in $a$.
\endproclaim

\demo{Proof of Claim 2.5}\enddemo

First note that by taking the saturation of $F$ in $E$ we
may assume that the quotient $E/F$ is also torsion free.

Take an $A$-slope-Harder-Narasimhan filtration of $F$
$$0 = F_0 \subset F_1 \subset \cdot\cdot\cdot \subset F_e
= F.$$
Set $Q_i = F_i/F_{i-1}$.  Then we have
$$h^0(F \otimes A^a) \leq \Sigma_i h^0(Q_i \otimes A^a),$$
and by the definition of an $A$-slope-Harder-Narasimhan
filtration we have for all $i = 1, 2,\cdot\cdot\cdot, e$
$$\mu_A(Q_i) = \frac{c_1(Q_i)}{rk(Q_i)} \cdot A \leq
\mu_A(Q_1) \leq \frac{c_1(E)}{rk(E)} \cdot A,$$
noting the second inequality comes from the fact that $E$
being $A$-Gieseker-\linebreak semistable implies $E$ being
$A$-slope-semistable.  We also have $$\Sigma_i rk(Q_i) =
rk(F).$$

Here we quote [Simpson92, Corollary1.7] without proof.

\proclaim{Lemma 2.6} There exists $b \in {\Bbb N}$
depending only on $r, X$ and $A$ such that for all
$A$-slope-semistable coherent torsion free sheaves
${\Cal F}$ of $rk({\Cal F}) \leq r$ on $X$
$$h^0({\Cal F}) \leq
\frac{rk({\Cal F})}{2}\text{deg}_AX\{\frac{\mu_A({\Cal
F})}{\text{deg}_AX} + b\}^2.$$
(Note also that
 if $\mu_A({\Cal F}) < 0$ then $h^0({\Cal F}) = 0$.)
\endproclaim

Going back to the proof of Claim 2.5, the lemma above
implies
$$\align
h^0(Q_i \otimes A^a) &\leq
\frac{rk(Q_i)}{2}\text{deg}_AX\{\frac{\mu_A(Q_i \otimes
A^a)}{\text{deg}_AX} + b\}^2\\
&=
\frac{rk(Q_i)}{2}\text{deg}_AX\{\frac{\mu_A(Q_i)}{\text{deg}_AX}
+ a + b\}^2.\\
\endalign$$
Thus we have
$$\align
h^0(F \otimes A^a) &\leq \Sigma
_ih^0(Q_i \otimes A^a)\\
&\leq \frac{rk(F) -
1}{2}\text{deg}_AX\{\frac{\mu_A(E)}{\text{deg}_AX} + a +
b\}^2\\
&+
\frac{1}{2}\text{deg}_AX\{\frac{\mu_A(Q_e)}{\text{deg}_AX}
+ a + b\}^2
\endalign$$
We remark that for all $\alpha \in {\Bbb Q}$ there exists
$\beta \in {\Bbb Q}$ and $a_{\beta} \in {\Bbb N}$ s.t. if
$$\frac{\mu_A(Q_e)}{\text{deg}_AX} \leq
\frac{\mu_A(E)}{\text{deg}_AX} - \beta$$
then for $a \geq a_{\beta}$
$$\align
&\frac{rk(F) -
1}{2}\text{deg}_AX\{\frac{\mu_A(E)}{\text{deg}_AX} + a +
b\}^2 +
\frac{1}{2}\text{deg}_AX\{\frac{\mu_A(Q_e)}{\text{deg}_AX}
+ a + b\}^2
\\
&\leq \frac{rk(F)}{2}\text{deg}_AX\{a - \alpha\}^2.
\endalign$$
For this we only have to take $\beta$ so that
$$\align
&\frac{rk(F) -
1}{2}\text{deg}_AX\{2(\frac{\mu_A(E)}{\text{deg}_AX} +
b)\} +
\frac{1}{2}\text{deg}_AX\{2(\frac{\mu_A(Q_e)}{\text{deg}_AX}
- \beta)\} \\
&\leq \frac{rk(F)}{2}\text{deg}_AX\{2(-\alpha)\}
\endalign$$
and then the existence of such $a_\beta$ is clear.

On the other hand, we can choose $\alpha$ and $
a_{\alpha} \in {\Bbb N}$ so that for all $a \geq
a_{\alpha}$ we have
$$\frac{1}{2}\text{deg}_AX(a -
\alpha)^2 < \frac{1}{rk(E)}\chi(E \otimes A^a).$$
Let $\beta$ be the number as above.  Then for $F \subset
E$ with
$$\frac{\mu_A(Q_e)}{\text{deg}_AX} \leq
\frac{\mu_A(E)}{\text{deg}_AX} - \beta$$
we have for $a \geq a_{\alpha}, a_{\beta}$
$$\frac{h^0(F \otimes A^a)}{rk(F)} < \frac{\chi(E
\otimes A^a)}{rk(E)}.$$
Moreover, since $S(r,c_1,c_2,A)$ is bounded, the set of
coherent sheaves $F$ such that $F \subset E$ for some $E
\in S(r,c_1,c_2,A)$ with $E/F$ being torsion free and that
$$\frac{\mu_A(Q_e)}{\text{deg}_AX} >
\frac{\mu_A(E)}{\text{deg}_AX} - \beta$$
is bounded (cf.Lemma 2.8).  Therefore, we may choose $a_{\gamma}$ so
that for all $a \geq a_{\gamma}$ we have
$$h^0(F \otimes A^a) = \chi(F \otimes A^a).$$
Furthermore, the set of the Hilbert polynomials for
all such $F$ is finite, so we may choose $a_{\delta}$
so that for all $a \geq a_{\delta}$ we have
$$\frac{\chi(F \otimes A^a)}{rk(F)} \leq \frac{\chi(E
\otimes A^a)}{rk(E)}$$
and that the equality holds iff
$$\frac{\chi(F \otimes A^a \otimes A^m)}{rk(F)} =
\frac{\chi(E \otimes A^a \otimes A^m)}{rk(E)} \text{\
for\ all\ }m \in {\Bbb Z},$$
i.e., iff
$$p(F,A,a) = p(E,A,a)$$
as polynomials in $a$.

Take $a_1 = \text{max} \{a_0, a_{\alpha}, a_{\beta},
a_{\gamma}, a_{\delta}\}$.

This completes the proof of Claim 2.5.

\vskip.1in

Now we go to the proof of the ``if" part of (ii) in the
Key GIT Lemma 2.4.

\vskip.1in

First we choose $a_1$ as in Claim 2.5 and take $a \geq
a_1$.

We take an ample line bundle $H \in Amp(X)_{\Bbb Q}$.

Suppose that a point
$$q:{\Cal O}_X^{\oplus l} \rightarrow (Univ)_q
\rightarrow 0$$
has the property that $(Univ)_q = E \otimes A^a$ is a
coherent torsion free sheaf of rank $rk(E) = r$, $c_1(E) =
c_1$, $c_2(E) = c_2$ such that $E$ is
$A$-Gieseker-semistable and that for all subsheaves $F
\subset E$ having the same averaged Euler characteristics
$p(F,A,a) = p(E,A,a)$ (as polynomials in $a$) we have
$$\frac{c_1(F)}{rk(F)} \cdot H \geq
\frac{c_1(E)}{rk(E)} \cdot H,$$
i.e.,
$$E \in S(r,c_1,c_2,A)_H,$$
and furthermore that
$$k^{\oplus l} = H^0({\Cal O}_X^{\oplus l}) \rightarrow
H^0((Univ)_q)$$
is an isomorphism.

By virtue of Lemma 2.2 we only have to show there exists
$M_H \in {\Bbb N}$ (bigger than the old $M_H$ in Lemma
2.2) such that for all $m \geq M_H$ we have
$$\frac{h^0(F \otimes A^a)}{\chi(F \otimes A^a \otimes
H^m)} \leq \frac{h^0(E \otimes A^a)}{\chi(E \otimes A^a
\otimes H^m)}$$
for all nonzero coherent subsheaves $F \subset E$ where
$F \otimes A^a$ is generated by some linear subspace of
$H^0(E \otimes A^a)$.

Once $a$ is fixed, the set of all such $F$ that $F
\subset E$ for some $E \in S(r,c_1,c_2,A)$ and that $F \otimes A^a$ is
generated by some linear subspace of
$H^0(E \otimes A^a)$ is bounded.  Thus the set of the
polynomials (in terms of $m$) $\chi(F \otimes A^a \otimes
H^m)$ is finite.  These polynomials all have the leading
term
$$\frac{rk(F)}{2}\text{deg}_HX \cdot m^2,$$
so we may choose $M_H$ big enough so that for $F$ with
$$\frac{h^0(F \otimes A^a)}{rk(F)} < \frac{h^0(E
\otimes A^a)}{rk(E)}$$
we get the desired inequality
$$\frac{h^0(F \otimes A^a)}{\chi(F \otimes A^a \otimes
H^m)} < \frac{h^0(E \otimes A^a)}{\chi(E \otimes A^a
\otimes H^m)}.$$
By the previous Claim 2.5, for those $F$ with the equality
$$\frac{h^0(F \otimes A^a)}{rk(F)} = \frac{h^0(E
\otimes A^a)}{rk(E)}$$
we have
$$\frac{\chi(F \otimes A^a \otimes A^m)}{rk(F)} =
\frac{\chi(E \otimes A^a \otimes A^m)}{rk(E)} \text{\
for\ all\ }m \in {\Bbb Z},$$
i.e.,
$$p(F,A,a) = p(E,A,a)$$
as polynomials in $a$.
By assumption, for those $F$ we have
$$\frac{c_1(F)}{rk(F)} \cdot H \geq \frac{c_1(E)}{rk(E)}
\cdot H,$$
and thus by rechoosing $M_H \in {\Bbb N}$ big enough if
necessary we have for all $m \geq M_H$
$$\align
\frac{\chi(F \otimes A^a \otimes H^m)}{rk(F)}
 &= \frac{1}{2}H^2m^2 + (\frac{c_1(F)}{rk(F)} \cdot H -
\frac{1}{2}K_X \cdot H)m + p(F,A,a)\\
&\geq \frac{1}{2}H^2m^2 + (\frac{c_1(E)}{rk(E)} \cdot H -
\frac{1}{2}K_X \cdot H)m + p(E,A,a)\\
&= \frac{\chi(E \otimes A^a \otimes H^m)}{rk(E)}\\
&= \frac{h^0(E \otimes A^a \otimes H^m)}{rk(E)} > 0.\\
\endalign$$
Thus again for those $F$ we have
$$\frac{h^0(F \otimes A^a)}{\chi(F \otimes A^a \otimes
H^m)} \leq \frac{h^0(E \otimes A^a)}{\chi(E \otimes A^a
\otimes H^m)}.$$
Therefore, the point $q \in Q$ is semistable with respect
to the action of $SL(l)$ and the linearization induced
from the Pl\"ucker embedding of $Grass(k^{\oplus l}
\otimes H^0(X,H^m),R_m)$.

This concludes the proof of the ``if" part.

\vskip.1in

Now we turn to the proof of the second half ``only if"
part of (ii).

Suppose that
$$q:{\Cal O}_X^{\oplus l} \rightarrow (Univ)_q = E
\otimes A^a \rightarrow 0$$
is a point of $Q$ which is semistable with respect to
the action of $SL(l)$ and the linearization induced
from the Pl\"ucker embedding of $Grass(k^{\oplus l}
\otimes H^0(X,H^m),R_m)$ for $m \geq M_H$ where $M_H$ is
the integer given in Lemma 2.3.

First we prove that $E$ is torsion free.

Let $T \subset E$ be the torsion of $E$.

We quote the following lemma of [Simpson92,Lemma1.17].

\proclaim{Lemma 2.7} Let
$$q:{\Cal O}_X^{\oplus l} \rightarrow (Univ)_q = E \otimes
A^a \rightarrow 0$$ be a point of $Q$ and $T \subset E$ be
the torsion part of $E$.  Then there exists a coherent
TORSION FREE sheaf $E'$ of rank $rk(E') = rk(E)$ such that
the Hilbert polynomials (in terms of $a$) are the same
$$\chi(E' \otimes A^a) = \chi(E \otimes A^a)$$
and that we have an inclusion
$$0 \rightarrow E/T \rightarrow E'.$$
\endproclaim

\vskip.1in

Take
$$0 \rightarrow E/T \rightarrow E'$$
as in Lemma 2.7.

Then for any quotient
$$E' \overset\psi\to\rightarrow G \rightarrow 0
\text{\ with\ }rk(G) > 0$$
we have
$$\align
\frac{h^0(G \otimes A^a)}{rk(G)} &\geq \frac{h^0(\psi(E/T)
\otimes A^a)}{rk(\psi(E/T))}\\
&\geq \frac{\chi(E \otimes A^a)}{rk(E)}\ (\text{by\
Lemma\ 2.3})\\
&= \frac{\chi(E' \otimes A^a)}{rk(E')}.\\
\endalign$$
Let $G$ be the $A$-slope-semistable quotient of
$E'$ with the smallest
$$\mu_A(G) = \frac{c_1(G)}{rk(G)} \cdot A,$$
in other words the last step quotient of an
$A$-slope-Harder-Narasimhan filtration of $E'$.

We apply Lemma 2.6 to conclude that there is $b \in {\Bbb
N}$ depending only on $r, X$ and $A$ such that
$$\align
\frac{\chi(E' \otimes A^a)}{rk(E')} &\leq \frac{h^0(G
\otimes A^a)}{rk(G)}\\
&\leq
\frac{1}{2}\text{deg}_AX\{\frac{\mu_A(G)}{\text{deg}_AX}
+ a + b\}^2.\\
\endalign$$
Since the polynomial in $a$
$$\chi(E' \otimes A^a) = \chi_{r,c_1,c_2}(A^a)$$
is fixed, there is a constant $c \in {\Bbb Q}$ and
$a_c\in {\Bbb N} (a_c \geq c)$ such that for all $a \geq
a_c$ we have $$\frac{\chi(E' \otimes A^a)}{rk(E')} \geq
\frac{1}{2}\text{deg}_AX(a - c)^2.$$
Thus
$$(a - c)^2 \leq \{\frac{\mu_A(G)}{\text{deg}_AX} + a +
b\}^2.$$
Since
$$\align
a - c \geq 0,\ \  &\frac{\mu_A(G
\otimes
A^a)}{\text{deg}_AX} = \frac{\mu_A(G)}{\text{deg}_AX} + a
\geq 0\\
&(\frac{h^0(G \otimes A^a)}{rk(G)} \geq
\frac{\chi_{r,c_1,c_2}(A^a)}{rk(E)} > 0)\\
\endalign$$
and since $b \in {\Bbb N}$, we have
$$a - c \leq \frac{\mu_A(G)}{\text{deg}_AX} + a + b,$$
which gives us
$$\mu_A(G) \geq \text{deg}_AX(- b - c).$$
Remark that $b$ and $c$ are constants taken independent
of $a$ or $m$.

\proclaim{Lemma 2.8} Let $L \in {\Bbb Q}$ be a fixed
 constant, $P(a)$ a fixed polynomial in $a$ and $A$ an
ample line bundle on $X$.  Then the set of all coherent
torsion free sheaves $E$ such that
$$\chi(E \otimes A^a) = P(a)$$
and that
$$\mu_A(G) \geq L$$
where $G$ is the last step quotient of an
$A$-slope-Harder-Narasimhan filtration of $E$, is bounded.
\endproclaim

\demo{Proof of Lemma 2.8}\enddemo

This is [Simpson92,Theorem1.1] proved using a lemma in [Maruyama81].
Here we give a simple alternative proof in the case where $X$ is a
surface, using Bogomolov inequality.

Take the Harder-Narasimhan filtration of $E$ with respect to
$A$-Gieseker-stability
$$0 = F_0 \subset F_1 \subset \cdot\cdot\cdot \subset
F_e = E.$$
(Note that the Harder-Narasimhan filtration of $E$ with respect to
$A$-slope-stability is obtained from the one above with respect to
$A$-Gieseker-stability by clustering together the factors with the same
averaged slope with respect to $A$.)

Set
$$G_i = F_i/F_{i-1}.$$
Write
$$c_1(G_i) \equiv \frac{c_1(G_i) \cdot A}{A^2}A  + {\Cal
L}_i$$
with
$${\Cal L} \in \frac{1}{A^2}\text{Pic}(X) \text{\ and\ } {\Cal L}_i \cdot A
= 0.$$
Then since
$$\mu_A(G_1) \geq \mu_A(G_2) \geq \cdot\cdot\cdot \mu_A(G_e) = \mu_A(G) \geq
L$$
is bounded from below and
$$rk(G_1)\mu_A(G_1) + rk(G_2)\mu_A(G_2) + \cdot\cdot\cdot + rk(G_e)\mu_A(G_e) =
c_1(E) \cdot A$$
is constant determined by the fixed polynomial $\chi(E \otimes A^a) = P(a)$,
there are only finite number of
 possibilities for $c_1(G_i) \cdot A$.

Now
$$\align
\chi(G_i \otimes A^a) &= \frac{1}{2}r_iA^2a^2 + \frac{1}{2}\{2c_1(G_i)\cdot
A - r_iK_X \cdot A\}a\\
&+ \frac{1}{2}\{c_1(G_i)^2 - 2c_2(G_i) - K_X \cdot c_1(G_i)\} +
r_i\chi({\Cal O}_X)\\
\endalign$$
and
$$\chi(E \otimes A^a) = \Sigma\chi(G_i \otimes A^a).$$
By looking at the constant term, we have
$$\align
&\text{the\ constant\ term\ of\ }P(a) - r\chi({\Cal O}_X)\\
&=\Sigma\frac{1}{2}\{c_1(G_i)^2 - 2c_2(G_i) - K_X\cdot c_1(G_i)\}\\
&\leq \Sigma\frac{1}{2}\{c_1(G_i)^2 - 2\frac{r_i-1}{2r_i}c_1(G_i)^2 -
K_X\cdot c_1(G_i)\} \text{\ by\ Bogomolov\ inequality}\\
&= \Sigma\frac{1}{2}\{\frac{1}{r_i}{\Cal L}_i^2 - K_X \cdot {\Cal L}_i +
\frac{1}{r_i}\frac{(c_1(G_i)\cdot A)^2}{A^2} - \frac{(c_1(G_i) \cdot
A)(K_X \cdot A)}{A^2}\}.\\
\endalign$$
Since we have finitely many possibilities for $c_1(G_i) \cdot A =
r_i\mu_A(G_i)$ and $r_i$, $K_X \cdot A$ and $A^2$ are constants, and since
${\Cal L}_i^2$ is quadratic in terms of ${\Cal L}_i \in A^{\perp}$ with
negative definite form, while $K_X \cdot {\Cal L}_i$ is linear, we conclude
$${\Cal L}_i \in A^{\perp} \cap \frac{1}{A^2}\text{Pic}(X) \subset
N^1(X)_{\Bbb Q}$$
has only a finite number of possibilities.  This implies $C_1(G_i)$ has
only a finite number of possibilities, and so does $c_2(G_i)$ again by
Bogomolov inequality.  Thus by [Gieseker77,Corollary1.3] the family of $G_i$
is bounded, and finally so is the family of $E$.

\vskip.1in

Applying this lemma, we conclude that the set of such
sheaves $E'$ remain in a bounded family which is
independent of $a$ or $m$.  Thus we may retake $a_1$ big
enough so that for all $a \geq a_1$ we have
$$h^0(E' \otimes A^a) = \chi(E' \otimes A^a)$$
and $E' \otimes A^a$ is generated by global sections.

Now taking the quotient sheaf in Lemma 2.3 to be $G =
E/T$, we get
$$h^0(E/T \otimes A^a) \geq \chi_{r,c_1,c_2}(A^a).$$
Since
$$\align
\chi(E' \otimes A^a) &= h^0(E' \otimes A^a)\\
&\geq h^0(E/T \otimes A^a)\\
&= \chi_{r,c_1,c_2}(A^a),\\
\endalign$$
we conclude that the natural inclusion
$$H^0(E/T \otimes A^a) \hookrightarrow H^0(E' \otimes
A^a)$$
is an isomorphism.  But since $E/T \hookrightarrow E'$
and since $H^0(E/T \otimes A^a) = H^0(E' \otimes
A^a)$ generates $E'$, we conclude $E/T = E'$.  Finally
since the Hilbert polynomials of $E$ and $E'$ are the
same, we conclude $T = 0$, i.e., $E$ is torsion free.
Furthermore, the equality
$$\align
h^0(E \otimes A^a) &= h^0(E/T \otimes A^a)\\
&= h^0(E' \otimes A^a) = \chi_{r,c_1,c_2}(A^a)\\
\endalign$$
together with the injectivity of Lemma 2.3
$$k^{\oplus l} \hookrightarrow H^0(E \otimes A^a)$$
implies
$$k^{\oplus l} \rightarrow H^0(E \otimes A^a)$$
is an isomorphism.

\vskip.1in

Next we prove that $E$ is $A$-Gieseker-semistable.
Suppose not.  Then there exists a nonzero torsion
free quotient $$E \rightarrow G \rightarrow 0$$
with $rk(G) > 0$ such that
$$\frac{\chi(G \otimes A^a)}{rk(G)} < \frac{\chi(E
\otimes A^a)}{rk(E)} \text{\ for\ }a >> 0.$$
Lemma 2.8 implies that the set of all such $G$ remains
in a bounded family independent of $a$ or $m$.
Therefore, we may choose $a_A \in {\Bbb N}$ ($a_A \geq
a_1,a_2$ where $a_2$ is the one in the second half of
Lemma 2.3) such that for all $a \geq a_A$ we have
$$\chi(G \otimes A^a) = h^0(G \otimes A^a).$$
Note also that the boundedness of all such $G$ implies
that the set of Hilbert polynomials $\chi(G \otimes
A^a \otimes H^m)$ as polynomials in $m$ is finite, and
that their leading terms are all
$$\frac{rk(G)}{2}
\text{deg}_HXm^2.$$
Thus we may choose $M_H$ (taken bigger than the one in
the second half of Lemma 2.3) such that for all $m \geq
M_H$ we have
$$\frac{h^0(G \otimes A^a)}{\chi(G \otimes A^a
\otimes H^m)} < \frac{\chi(E \otimes A^a)}{\chi(E
\otimes A^a \otimes H^m)},$$
contradicting the conclusion of the second half of
Lemma 2.3.

Finally let $F \subset E$ be a subsheaf of $E$ having the
same averaged Euler characteristics
$$p(F,A,a) = p(E,A,a) (\text{\ as\ polynomials\ in\ }a).$$
Suppose
$$\frac{c_1(F)}{rk(F)} \cdot H < \frac{c_1(E)}{rk(E)}
\cdot H.$$
Let $G$ be the quotient
$$0 \rightarrow F \rightarrow E \rightarrow G \rightarrow
0.$$
Note that the conditions imply that $G$ is also a
coherent torsion free sheaf which is
$A$-Gieseker-semistable with
$$p(G,A,a) = p(E,A,a) (\text{\ as\ polynomials\ in\
}a),$$
and that
$$\frac{c_1(G)}{rk(G)} \cdot H > \frac{c_1(E)}{rk(E)}
\cdot H.$$
Then for $m \in {\Bbb N}$ we have
$$\align
p(G \otimes A^a,H,m) &= \frac{1}{2}H^2m^2 +
(\frac{c_1(G)}{rk(G)} \cdot H - \frac{1}{2}K_X \cdot H)m
+ p(G,A,a)\\
&> \frac{1}{2}H^2m^2 +
(\frac{c_1(E)}{rk(E)} \cdot H - \frac{1}{2}K_X \cdot H)m
+ p(E,A,a)\\
&= p(E \otimes A^a,H,m)\\
\endalign$$
Since the set of all such $G$ is bounded independent of
$a$, we may choose $a_A$ so that for all $a \geq a_A$

$$\frac{h^0(G \otimes A^a)}{rk(G)} = p(G,A,a) = p(E,A,a) =
\frac{h^0(E \otimes A^a)}{rk(E)} > 0.$$
Again since the set of all such $G$ is bounded, we may
assume that for $m \geq M_H$ we have
$$ 0 < \frac{h^0(G \otimes A^a \otimes H^m)}{rk(G)} = p(G
\otimes A^a,H,m) > p(E \otimes A^a,H,m).$$
These imply
$$\align
&\frac{h^0(G \otimes A^a)}{\chi(G \otimes A^a \otimes
H^m)} = \frac{p(G,A,a)}{p(G \otimes A^a,H,m)}\\
&< \frac{p(E,A,a)}{p(E \otimes A^a,H,m)} = \frac{\chi(E
\otimes A^a)}{\chi(E \otimes A^a \otimes H^m)},\\
\endalign$$
contradicting again the conclusion of the second half of
Lemma 2.3.

This completes the proof of the Key GIT Lemma.

\vskip.1in

\S 3. Rationally Twisted Gieseker-Semistability and Stratifications

\vskip.1in

In this section, we analyse the change of the set of
the $A$-Gieseker-semistable sheaves (the polarization $A$ being fixed)
when tensored by a line bundle ${\Cal L} \in \text{Pic}(X)$.  The
analysis naturally leads us to the notion of Gieseker-semistability
twisted by a rational line bundle ${\Cal L} \in \text{Pic} \otimes
{\Bbb Q}$.  We will give a stratification of the space of twists,
which describes the change of the set of ${\Cal L}$-twisted
$A$-Gieseker-semistable sheaves as ${\Cal L}$ varies, based upon
the stratification of $\Delta - \{0\} = \coprod_s\Delta_s$.
  These two stratifications not only completely describe the change of
$S(r,c_1,c_2,H)$ as $H$ varies, but also illustrate how the
transformation is factorized into a sequence of flips. \vskip.1in

Consider the set of all coherent torsion free sheaves $F$ on $X$ as
in Proposition 1.6 such that there exists an exact sequence
$$0 \rightarrow F \rightarrow E \rightarrow G \rightarrow 0$$
where $E \in \mu(r,c_1,c_2,H)$ for some $H \in \Delta - \{0\}$, and
where the quotient $G$ is a coherent torsion free sheaf, satisfying
$$\frac{c_1(F)}{rk(F)} \cdot H = \frac{c_1(E)}{rk(E)} \cdot H$$
but
$$0 \neq \frac{c_1(F)}{rk(F)} - \frac{c_1(E)}{rk(E)} \text{\ in\
}N^1(X)_{\Bbb Q}.$$

Take the set of all hyperplanes
$$L_F := \{z \in N^1(X)_{\Bbb Q};z \cdot (\frac{c_1(F)}{rk(F)} -
\frac{c_1(E)}{rk(E)}) = 0\}$$
for all such sheaves $F$.  Proposition 1.6 implies that $\{L_F\}$
is a finite set.

\proclaim{Definition 3.1} A $d$-cell of $\Delta - \{0\}$ associated
to $\{L_F\}$ ($d = 1, 2, \cdot\cdot\cdot, \rho_{\Delta} =
\text{dim}_{\Bbb Q}V(\Delta)$, where $V(\Delta)$ is the vector
subspace of $N^1(X)_{\Bbb Q}$ generated by $\Delta$.) is a
connected component of
$$\Delta \cap L_{F_1} \cap \cdot\cdot\cdot L_{F_{\rho_{\Delta} -
d}} - (\cup_{F \neq F_i, i = 1, \cdot\cdot\cdot, \rho_{\Delta} -
d}L_F).$$
\endproclaim

Now we have a stratification of $\Delta - \{0\}$
$$\Delta - \{0\} = \coprod \Delta_s$$
into the $d$-cells where $d = 1, \cdot\cdot\cdot, \rho_{\Delta}$.

Before going into further discussion, we introduce the notion of
rationally twisted Gieseker-semistability, which just slightly
generalizes the classical notion of \hfil\linebreak
Gieseker-semistability but
provides us with the right category to factorize the transformations
among the various moduli spaces.

\proclaim{Definition 3.2} Let $H$ be an ample line bundle.  A coherent
torsion free sheaf $E$ is said to be ${\Cal L}$-twisted
$H$-Gieseker-semistable (resp. ${\Cal L}$-twisted
$H$-Gieseker-stable) for ${\Cal L} \in \text{Pic}(X) \otimes {\Bbb Q}$
iff for all $F \subset E$ (resp. for all $F
\underset{\neq}\to\subset E$)
$$\align
\frac{\chi(F \otimes {\Cal L} \otimes
H^n)}{rk(F)} &\leq \frac{\chi(E \otimes {\Cal L} \otimes
H^n)}{rk(E)} \text{\ for\ }n >> 0\\
&(\text{resp.\ }<)\\
\endalign
$$
where we compute the Euler characteristics formally using the
Riemann-Roch formula, e.g.,
$$\align
\frac{\chi(E \otimes {\Cal L} \otimes H^n)}{rk(E)} &=
\frac{1}{2}H^2n^2\\
&+ \{\frac{c_1(E)}{rk(E)} \cdot H + {\Cal L} \cdot H -
\frac{1}{2}K_X \cdot H\}n\\
&+ \frac{1}{2}{\Cal L}^2 + \frac{c_1(E)}{rk(E)} \cdot {\Cal L} -
\frac{1}{2}K_X \cdot {\Cal L}\\
&+ \frac{1}{2}\frac{c_1(E)^2 - 2c_2(E) - c_1(E) \cdot K_X}{rk(E)}\\
&+ \chi({\Cal O}_X).\\
\endalign$$
We denote by $S((r,c_1,c_2) \otimes {\Cal L},H)$ the set of all
coherent torsion free sheaves $E$ of $rk(E) = r, c_1(E) = c_1,
c_2(E) = c_2$, that are ${\Cal L}$-twisted $H$-Gieseker-semistable.
\endproclaim

\proclaim{Remark 3.3}\endproclaim

\ \ (i) As is clear from the definition, ${\Cal L}$-twistedness only
depends on the class of ${\Cal L}$ in $N^1(X)_{\Bbb Q}$.

\ \ (ii) When ${\Cal L}$ is an INTEGRAL line bundle ${\Cal L} \in
\text{Pic}(X)$, then $S((r,c_1,c_2) \otimes {\Cal L},H)$ can be
identified with $S(r,{c_1}_{\Cal L},{c_2}_{\Cal L},H)$ where
$$\align
{c_1}_{\Cal L} &= c_1 + rc_1({\Cal L})\\
{c_2}_{\Cal L} &= c_2 + (r - 1)c_1 \cdot c_1({\Cal L}) + \frac{r(r -
1)}{2}c_1({\Cal L})^2\\
\endalign$$
under the one to one correspondence
$$E \in S((r,c_1,c_2) \otimes {\Cal L},H) \leftrightarrow E \otimes
{\Cal L} \in S(r,{c_1}_{\Cal L},{c_2}_{\Cal L},H).$$

\ \ (iii) The Harder-Narasimhan filtration with respect to
rationally twisted \linebreak
Gieseker-(semi)stability can be taken just as in
the classical case and the notion of ${\Cal L}$-twisted Seshadri
equivalence is also well defined.

The stratification $\Delta - \{0\} = \coprod \Delta_s$ gives a basis
to study the change of \linebreak
$S(r,c_1,c_2,H)$ when $H$ varies and the change
of $S((r,c_1,c_2) \otimes {\Cal L},H)$ when ${\Cal L}$ changes, as is
indicated by the next lemma, whose proof is immediate from the
construction.

\proclaim{Lemma 3.4}

(i) For $H, H' \in \Delta_s$, we have
$$\align
\mu(r,c_1,c_2,H) &= \mu(r,c_1,c_2,H')\\
S(r,c_1,c_2,H) &= S(r,c_1,c_2,H').\\
\endalign$$
Moreover, the Seshadri equivalence classes of $S(r,c_1,c_2,H)$ with
respect to $H$-Gieseker-stability are the same as the Seshadri
equivalence classes of $S(r,c_1,c_2,H')$ with respect to
$H'$-Gieseker-stability.

(ii) Fix $H \in \Delta_s$.  Then for any ${\Cal L} \in
V(\Delta_s)$
, we have
$$S(r,c_1,c_2,H) = S((r,c_1,c_2) \otimes {\Cal L},H).$$
Moreover, the Seshadri equivalence classes of $S(r,c_1,c_2,H)$ with
respect to $H$-Gieseker-stability are the same as the Seshadri
equivalence classes of $S((r,c_1,c_2) \otimes {\Cal L},H)$ with
respect to ${\Cal L}$-twisted $H$-Gieseker-stability.
\endproclaim

Let $\Delta_s$ and $\Delta_{s'}$ be two $d$-cells such that
$V(\Delta_s) = V(\Delta_{s'})$ and that $\Delta_s$ and $\Delta_{s'}$
 are separated by a $d - 1$-cell $W$, i.e., $W$ is the unique $d -
1$-cell contained in $\overline{\Delta_s} \cap
\overline{\Delta_{s'}}$.  We take an ample line bundle $A \in W$.

\proclaim{Proposition 3.5} There exists a stratification
$$V(\Delta_s) = \coprod L_i \coprod M_j$$
consisting of a finite number of hyperplanes $L_i$ parallel to
$V(W)$ and connected components $M_j$ of $V(\Delta_s) - (\coprod
L_i)$, which determines the change of $S((r,c_1,c_2) \otimes {\Cal
L},A)$ as ${\Cal L}$ varies in $V(\Delta_s)$, i.e., for ${\Cal L},
{\Cal L}' \in V(\Delta_s)$,
$$S((r,c_1,c_2) \otimes {\Cal L},A) = S((r,c_1,c_2) \otimes {\Cal
L}',A)$$
if and only if ${\Cal L}$ and ${\Cal L}'$ belong to the same
stratum $L_i$ or $M_j$.  This stratification is independent of the
choice of $A \in W$.

Moreover,
$$S(r,c_1,c_2,A) = S(r,c_1,c_2,H) \text{\ for\ }H \in
\Delta_s$$
if and only if $\overline{\Delta_s}$ is contained in one of the
strata $M_j$.  By induction on the dimension $d$ of the $d$-cells,
these rules completely determine the change of the set
$S(r,c_1,c_2,H)$ as $H$ varies in $\Delta - \{0\}$.
\endproclaim

\demo{Proof of Proposition 3.5}\enddemo

The condition for a coherent torsion free sheaf
$$E \in \mu(r,c_1,c_2,A)$$
to be ${\Cal L}$-twisted $A$-Gieseker-semistable is that
$$(\diamondsuit) \frac{\chi(F \otimes {\Cal L} \otimes A^n)}{rk(F)}
\leq \frac{\chi(E \otimes {\Cal L} \otimes A^n)}{rk(E)} \text{\ for\
}n >> 0$$
for all nonzero subsheaves $F \subset E$.  By taking the saturation
we only have to check the condition for those subsheaves $F$ with
the torsion free quotient $E/F$.

If there is a subsheaf $F \subset E$ with
$$\{\frac{c_1(F)}{rk(F)} - \frac{c_1(E)}{rk(E)}\} \cdot A > 0,$$
then the condition $(\diamondsuit)$ is never satisfied for any
${\Cal L} \in V(\Delta_s)$ for this $F$ and $E$ is not ${\Cal
L}$-twisted $A$-Gieseker-semistable for any ${\Cal L} \in
V(\Delta_s)$.

For a subsheaf $F \subset E$ with
$$\{\frac{c_1(F)}{rk(F)} - \frac{c_1(E)}{rk(E)}\} \cdot A < 0,$$
the condition $(\diamondsuit)$ is always satisfied for any ${\Cal L} \in
V(\Delta_s)$.

Therefore, the only subsheaves $F \subset E$ for which the validity
of the condition $(\diamondsuit)$ depends on ${\Cal L}$ are the ones
with
$$\{\frac{c_1(F)}{rk(F)} - \frac{c_1(E)}{rk(E)}\} \cdot A = 0.$$
For such $F$, the condition $(\diamondsuit)$ is equivalent to
$$\align
&\{\frac{c_1(F)}{rk(F)} - \frac{c_1(E)}{rk(E)}\} \cdot {\Cal L}\\
&+ \frac{1}{2}\frac{c_1(F)^2 - 2c_2(F) - c_1(F) \cdot K_X}{rk(F)}\\
&- \frac{1}{2}\frac{c_1(E)^2 - 2c_2(E) - c_1(E) \cdot K_X}{rk(E)}\\
&\leq 0.\\
\endalign$$
If $\{\frac{c_1(F)}{rk(F)} - \frac{c_1(E)}{rk(E)}\}$ is numerically
equivalent to 0, then again the condition $(\diamondsuit)$ is
independent of ${\Cal L} \in V(\Delta_s)$.

Take the set of hyperplanes
$$\align
M_F = \{z \in V(\Delta_s);&(\frac{c_1(F)}{rk(F)} -
\frac{c_1(E)}{rk(E)}) \cdot z\\
&+ \frac{1}{2}\frac{c_1(F)^2 - 2c_2(F) - c_1(F) \cdot K_X}{rk(F)}\\
&- \frac{1}{2}\frac{c_1(E)^2 - 2c_2(E) - c_1(E) \cdot K_X}{rk(E)} =
0\}\\
\endalign$$
for all subsheaves $F \subset E$ for some $E \in \mu(r,c_1,c_2,A)$
with $E/F$ torsion free,
$$\{\frac{c_1(F)}{rk(F)} - \frac{c_1(E)}{rk(E)}\} \cdot A = 0$$
and $\{\frac{c_1(F)}{rk(F)} - \frac{c_1(E)}{rk(E)}\}$ being not
numerically equivalent to 0.  Proposition 1.6 says that $\{M_F\}$ is
a finite set.  Note that the $M_F$ are all parallel to $V(W)$, since
$A \in W$ where $W$ is one of the $d - 1$-cells of the stratification
of $\Delta - \{0\}$.

Fix $E \in \mu(r,c_1,c_2,A)$.  Then the previous argument shows
that the locus of ${\Cal L}$ in $V(\Delta_s)$ where $E$ is ${\Cal
L}$-twisted $A$-Gieseker-semistable is either empty, one of the
hyperplanes $M_F$, a closed subspace sandwiched between two of
the hyperplanes $M_F$, a closed half subspace whose boundary is
one of the hyperplanes $M_F$ or the entire space $V(\Delta_s)$.

Let $L$ be one of the boundary hyperplanes and ${\Cal L} \in L$.
Since ${\Cal L}$ is on the boundary, there exists a nonzero
proper subsheaf $F \subset E$ with torsion free $E/F$ such that
$$\frac{\chi(F \otimes {\Cal L} \otimes A^n)}{rk(F)} = \frac{\chi(E
\otimes {\Cal L} \otimes A^n)}{rk(E)} = \frac{\chi(E/F \otimes {\Cal
L} \otimes A^n)}{rk(E/F)} \text{\ for\ all\ }n,$$
that both $F$ and $E/F$ are ${\Cal L}$-twisted
$A$-Gieseker-semistable and that
$$\{\frac{c_1(F)}{rk(F)} - \frac{c_1(E)}{rk(E)}\} \cdot {\Cal M}
\overset<\to{\underset>\to=} 0$$
if and only if
$$\{\frac{c_1(E/F)}{rk(E/F)} - \frac{c_1(E)}{rk(E)}\} \cdot {\Cal M}
\overset>\to{\underset<\to=} 0.$$
with
$$\{\frac{c_1(F)}{rk(F)} - \frac{c_1(E)}{rk(E)}\} (\text{\ and\ thus\
also\ }\{\frac{c_1(E/F)}{rk(E/F)} - \frac{c_1(E)}{rk(E)}\})$$
not being numerically trivial on $V(\Delta)$.  This implies that the
sheaf $$F \oplus E/F$$ is ${\Cal M}$-twisted $A$-Gieseker-semistable
if and only if ${\Cal M} \in L$.

This completes the argument for the proof that there is a
stratification as desired $$V(\Delta_s) = \coprod L_i \coprod M_j$$
consisting of a finite number of hyperplanes $L_i$ (a part of the
$M_F$) parallel to $V(W)$ and connected components $M_j$ of
$V(\Delta_s) - (\coprod L_i)$ which determines the change of the
set $S((r,c_1,c_2) \otimes {\Cal L},A)$ as ${\Cal L}$ varies in
$V(\Delta_s)$.

The last part of the proposition follows immediately from the
fact that for any $H \in \Delta_s$
$$S(r,c_1,c_2,H) = S((r,c_1,c_2) \otimes H^n,A) \text{\ for\ }n >>
0,$$
which will be shown in the next lemma.

This completes the proof of Proposition 3.5.

\vskip.1in

We remark that since for any ${\Cal L} \in V(\Delta_s)$ we have
$$S((r,c_1,c_2) \otimes {\Cal L},A) = S((r,c_1,c_2) \otimes {\Cal L}
\otimes A^a,A)$$
and since ${\Cal L} \otimes A^a$ is ample for $a >> 0$, each stratum
of the stratification above has an ample representative, and that it is
actually determined by its intersection with the ample cone
$$V(\Delta_s) \cap Amp(X)_{\Bbb Q} = (\coprod L_i \coprod M_j) \cap
Amp(X)_{\Bbb Q}.$$

\proclaim{Lemma 3.6} Take an ample line bundle $H \in \Delta_s$
(resp. $H' \in \Delta_{s'}$).  Then
$$S((r,c_1,c_2) \otimes H^n,A) = S(r,c_1,c_2,H) \text{\ for\ }n >>
0$$
$$(\text{resp.\ } S((r,c_1,c_2) \otimes {H'}^{n'},A) = S(r,c_1,c_2,H')
\text{\ for\ }n' >> 0).$$
Moreover, the Seshadri equivalence classes of $S((r,c_1,c_2) \otimes
H^n,A)$ (resp. \hfil\linebreak $S((r,c_1,c_2) \otimes
{H'}^{n'},A)$) with respect to $H^n$-twisted (resp.
${H'}^{n'}$-twisted ) $A$-Gieseker-stability are the same as the
Seshadri equivalence classes of $S(r,c_1,c_2,H)$ (resp.
$S(r,c_1,c_2,H')$) with respect to $H$-Gieseker-stability (resp.
$H'$-Gieseker-stability).
\endproclaim

\demo{Proof of Lemma 3.6}\enddemo

Let $P = \{(c_1(F),c_2(F))\}$ be the set of the pairs of numerical
classes for such sheaves $F \subset E$ with $E/F$ torsion free for
some $E \in \mu(r,c_1,c_2,A)$ that
$$\{\frac{c_1(F)}{rk(F)} - \frac{c_1(E)}{rk(E)}\}\cdot A = 0.$$
$P$ is a finite set by Proposition 1.6, and thus we can choose $n >>
0$ so that for all $(c_1(F),c_2(F)) \in P$, if
$$\{\frac{c_1(F)}{rk(F)} - \frac{c_1(E)}{rk(E)}\} \cdot H < 0$$
then
$$\align
&\{\frac{c_1(F)}{rk(F)} - \frac{c_1(E)}{rk(E)}\} \cdot H^n\\
&+ \frac{1}{2}\frac{c_1(F)^2 - 2c_2(F) - c_1(F) \cdot K_X}{rk(F)}\\
&- \frac{1}{2}\frac{c_1(E)^2 - 2c_2(E) - c_1(E) \cdot K_X}{rk(E)} <
0,\\
\endalign$$
and if
$$\{\frac{c_1(F)}{rk(F)} - \frac{c_1(E)}{rk(E)}\} \cdot H > 0$$
then
$$\align
&\{\frac{c_1(F)}{rk(F)} - \frac{c_1(E)}{rk(E)}\} \cdot H^n\\
&+ \frac{1}{2}\frac{c_1(F)^2 - 2c_2(F) - c_1(F) \cdot K_X}{rk(F)}\\
&- \frac{1}{2}\frac{c_1(E)^2 - 2c_2(E) - c_1(E) \cdot K_X}{rk(E)}
> 0.\\ \endalign$$
Take $E \in S(r,c_1,c_2,H)$.  Then for all subsheaves $F \subset E$
($E/F$ torsion free) we have either
$$\{\frac{c_1(F)}{rk(F)} - \frac{c_1(E)}{rk(E)}\} \cdot H < 0$$
or
$$\{\frac{c_1(F)}{rk(F)} - \frac{c_1(E)}{rk(E)}\} \cdot H = 0$$
together with the condition
$$\align
&\frac{1}{2}\frac{c_1(F)^2 - 2c_2(F) - c_1(F) \cdot K_X}{rk(F)}\\
&- \frac{1}{2}\frac{c_1(E)^2 - 2c_2(E) - c_1(E) \cdot K_X}{rk(E)}
\leq 0.\\
\endalign$$
In the first case, we have
$$\{\frac{c_1(F)}{rk(F)} - \frac{c_1(E)}{rk(E)}\} \cdot A \leq 0.$$
If the strict inequality $<$ holds, then the condition
$(\diamondsuit)$ is immediate.  If the equality $=$ holds, then the
choice of $n$ guarantees $(\diamondsuit)$.

In the second case, we have
$$\{\frac{c_1(F)}{rk(F)} - \frac{c_1(E)}{rk(E)}\} \cdot A = 0.$$
Since
$$\align
&\frac{1}{2}\frac{c_1(F)^2 - 2c_2(F) - c_1(F) \cdot K_X}{rk(F)}\\
&- \frac{1}{2}\frac{c_1(E)^2 - 2c_2(E) - c_1(E) \cdot K_X}{rk(E)}
\leq 0,\\
\endalign$$
again the condition $(\diamondsuit)$ is satisfied.

Therefore, we conclude $E \in S((r,c_1,c_2) \otimes H^n,A)$.

Now suppose $E \in S((r,c_1,c_2) \otimes H^n,A)$.  Note that $E \in
\mu(r,c_1,c_2,A)$.  If $E \notin S(r,c_1,c_2,H)$, then there exists
a subsheaf $F \subset E$ such that either
$$\{\frac{c_1(F)}{rk(F)} - \frac{c_1(E)}{rk(E)}\} \cdot H > 0$$
or
$$\{\frac{c_1(F)}{rk(F)} - \frac{c_1(E)}{rk(E)}\} \cdot H = 0$$
together with the condition
$$\align
&\frac{1}{2}\frac{c_1(F)^2 - 2c_2(F) - c_1(F) \cdot K_X}{rk(F)}\\
&- \frac{1}{2}\frac{c_1(E)^2 - 2c_2(E) - c_1(E) \cdot K_X}{rk(E)}
> 0.\\
\endalign$$

In the first case, if
$$\{\frac{c_1(F)}{rk(F)} - \frac{c_1(E)}{rk(E)}\} \cdot A > 0,$$
then $E \notin S((r,c_1,c_2) \otimes {H^n},A)$, a contradiction!  If
$$\{\frac{c_1(F)}{rk(F)} - \frac{c_1(E)}{rk(E)}\} \cdot A = 0,$$
then the choice of $n$ implies
$$\align
&\{\frac{c_1(F)}{rk(F)} - \frac{c_1(E)}{rk(E)}\} \cdot H^n\\
&+ \frac{1}{2}\frac{c_1(F)^2 - 2c_2(F) - c_1(F) \cdot K_X}{rk(F)}\\
&- \frac{1}{2}\frac{c_1(E)^2 - 2c_2(E) - c_1(E) \cdot K_X}{rk(E)}
> 0,\\
\endalign$$
and thus $E \notin S((r,c_1,c_2) \otimes {H^n},A)$, again a
contradiction!

In the second case, we have
$$\{\frac{c_1(F)}{rk(F)} - \frac{c_1(E)}{rk(E)}\} \cdot A = 0,$$
and
$$\align
&\frac{1}{2}\frac{c_1(F)^2 - 2c_2(F) - c_1(F) \cdot K_X}{rk(F)}\\
&- \frac{1}{2}\frac{c_1(E)^2 - 2c_2(E) - c_1(E) \cdot K_X}{rk(E)}
> 0\\
\endalign$$
implies
$$\align
&\{\frac{c_1(F)}{rk(F)} - \frac{c_1(E)}{rk(E)}\} \cdot H^n\\
&+ \frac{1}{2}\frac{c_1(F)^2 - 2c_2(F) - c_1(F) \cdot K_X}{rk(F)}\\
&- \frac{1}{2}\frac{c_1(E)^2 - 2c_2(E) - c_1(E) \cdot K_X}{rk(E)}
> 0.\\
\endalign$$
Thus $E \notin S((r,c_1,c_2) \otimes {H^n},A)$, a contradiction!

Therefore, we conclude $E \in S(r,c_1,c_2,H)$.

The statement about the Seshadri equivalence classes follows
similarly.

\vskip.2in

Let $M_0$ and $M_{l+1}$ be the two strata containing $H^n$ (for $n
>> 0$) and ${H'}^{n'}$ (for $n' >> 0$) respectively.  Then Lemma 3.6
indicates that at least set-theoretically these two strata
correspond to the moduli spaces $M(r,c_1,c_2,H)$,
$M(r,c_1,c_2,H')$.  One would expect that going from $M_0$ in the
stratification $V(\Delta_s) = \coprod L_i \coprod M_j$ to $M_{l+1}$
corresponds to a sequence of moduli spaces of rationally twisted
$A$-Gieseker-semistable sheaves.  It should be noted that some of
the strata may not contain any integral points corresponding to
$\text{Pic}(X)$ and that this is the reason we are naturally led to
introduce the notion of rationally twisted Gieseker-semistability.

\S 4 and \S 5 will be
devoted to showing not only set-theoretically but
scheme-theoretically that this is indeed the case.  We construct the
moduli space of rationally twisted Gieseker-semistable
sheaves.  Moreover we show that there is a sequence of flips among
them, each of which is governed by the Mumford-Thaddeus principle of
GIT: Let $M_0, L_0, L_1, M_1, \cdot\cdot\cdot, L_l, M_{l+1}$ be a
sequence of strata starting with $M_0$ containing $H^n$ (for $n >>
0$) and ending with $M_{l+1}$ containing ${H'}^{n'}$ (for $n' >>
0)$ such that
$$\overline{M_i} \cap \overline{M_{i+1}} = L_i$$
for $i = 0, 1, \cdot\cdot\cdot, l$.  (Note that these strata
actually exhaust all the strata in $V(\Delta_s) = \coprod L_i
\coprod M_j$.)  Choose representatives ${\Cal L}_i, {\Cal M}_i \in
\text{Pic}(X)_{\Bbb Q}$ with ${\Cal L}_i \in L_i, {\Cal M}_i \in
M_i$.

Our aim is to show that there is a sequence of Thaddeus-type
flips
$$
\CD
M((r,c_1,c_2) \otimes {\Cal M}_i,A) @.@.@.@.M((r,c_1,c_2)
\otimes {\Cal M}_{i+1},A)\\
@.\searrow @.@.\swarrow @.\\
@.@.M((r,c_1,c_2) \otimes {\Cal L}_i,A)@.@.\\
\endCD
$$
for $i = 1, 2, \cdot\cdot\cdot, l$, each of which is
a transformation constructed by the Key GIT Lemma of \S 2
and thus governed by the Mumford-Thaddeus principle.

\vskip.1in

\S 4. Construction of Flip: Integral Case

\vskip.1in

We use the notation in \S 3.

In this section, we construct the flip
$$
\CD
M((r,c_1,c_2) \otimes {\Cal M}_i,A) @.@.@.@.M((r,c_1,c_2)
\otimes {\Cal M}_{i+1},A)\\
@.\searrow @.@.\swarrow @.\\
@.@.M((r,c_1,c_2) \otimes {\Cal L}_i,A)@.@.\\
\endCD
$$
when a representative ${\Cal L}_i \in L_i$ can be taken to be an INTEGRAL line
bundle, i.e., ${\Cal L} \in
\text{Pic}(X)$, using the Key GIT Lemma of \S 2.  (We will deal with
the more delicate case of the construction of the flip when ${\Cal
L}_i \in L_i$ is only a RATIONAL line bundle in the next section.)

\proclaim{Theorem 4.1}

(i) The one to one correspondence
$$E \in S((r,c_1,c_2) \otimes {\Cal L}_i,A) \leftrightarrow E
\otimes {\Cal L}_i \in S(r,{c_1}_{{\Cal L}_i},{c_2}_{{\Cal L}_i},A)$$
gives the one to one correspondence between
$$S((r,c_1,c_2) \otimes {\Cal M}_i,A) \text{\ and\ }S(r,{c_1}_{{\Cal
L}_i},{c_2}_{{\Cal L}_i},A)_{H'}$$
and that between
$$S((r,c_1,c_2) \otimes {\Cal M}_{i+1},A) \text{\ and\
}S(r,{c_1}_{{\Cal L}_i},{c_2}_{{\Cal L}_i},A)_H.$$

(ii) $M(r,{c_1}_{{\Cal L}_i},{c_2}_{{\Cal L}_i},A)_{H'}$ (resp.
$M(r,{c_1}_{{\Cal L}_i},{c_2}_{{\Cal
L}_i},A)$, $M(r,{c_1}_{{\Cal L}_i},{c_2}_{{\Cal L}_i},A)_H$)
constructed in \S 2 gives the moduli space $M((r,c_1,c_2) \otimes
{\Cal M}_i,A)$ (resp. $M((r,c_1,c_2) \otimes
{\Cal L}_i,A)$, $M((r,c_1,c_2) \otimes
{\Cal M}_{i+1},A)$) which coarsely represents the Seshadri
equivalence classes of $S((r,c_1,c_2) \otimes {\Cal M}_i,A)$ (resp.
$S((r,c_1,c_2) \otimes {\Cal L}_i,A)$, $S((r,c_1,c_2) \otimes {\Cal
M}_{i+1},A)$) with respect to ${\Cal M}_i$-twisted (resp.
${\Cal L}_i$-twisted, ${\Cal M}_{i+1}$-twisted)
$A$-Gieseker-stability.

(iii) The diagram of morphisms
$$
\CD
M({c_1}_{{\Cal L}_i},{c_2}_{{\Cal L}_i},A)_{H'} @.@.@.@.M((r,{c_1}_{{\Cal
L}_i},{c_2}_{{\Cal
L}_i},A)_H\\ @.\searrow @.@.\swarrow @.\\
@.@.M(r,{c_1}_{{\Cal L}_i},{c_2}_{{\Cal L}_i},A)@.@.\\
\endCD
$$
which arises from the Key GIT Lemma gives the desired flip
$$
\CD
M((r,c_1,c_2) \otimes {\Cal M}_i,A) @.@.@.@.M((r,c_1,c_2)
\otimes {\Cal M}_{i+1},A)\\
@.\psi_i\searrow @.@.\swarrow\psi_i^+ @.\\
@.@.M((r,c_1,c_2) \otimes {\Cal L}_i,A)@.@.\\
\endCD
$$
\endproclaim

\demo{Proof of Theorem 4.1}\enddemo

(i) Suppose $E \in S((r,c_1,c_2) \otimes {\Cal M}_i,A)$.  Since
${\Cal L}_i$ is on the boundary $L_i$ of the stratum $M_i$
containing ${\Cal M}_i$,
$$S((r,c_1,c_2) \otimes {\Cal M}_i,A) \subset S((r,c_1,c_2) \otimes
{\Cal L}_i,A).$$
Therefore, $E \otimes {\Cal L}_i$ is $A$-Gieseker-semistable.  If $F
\otimes {\Cal L}_i$ is a subsheaf of $E \otimes {\Cal L}_i$ with
$$p(F \otimes {\Cal L}_i,A,a) = \frac{\chi(F \otimes {\Cal L}_i
\otimes A^a)}{rk(F)} = \frac{\chi(E \otimes {\Cal L}_i \otimes
A^a)}{rk(E)} = p(E \otimes {\Cal L}_i,A,a)$$
and
$$\frac{c_1(F)}{rk(F)} \cdot H' < \frac{c_1(E)}{rk(E)} \cdot H',$$
then
$$\frac{c_1(F)}{rk(F)} \cdot A = \frac{c_1(E)}{rk(E)} \cdot A$$
$$\frac{1}{2}\frac{c_1(F)^2 - 2c_2(F) - c_1(F) \cdot K_X}{rk(F)} =
\frac{1}{2}\frac{c_1(E)^2 - 2c_2(E) - c_1(E) \cdot K_X}{rk(E)}$$
and
$$\frac{c_1(F)}{rk(F)} \cdot {\Cal M}_i > \frac{c_1(E)}{rk(E)}
\cdot {\Cal M}_i.$$
This implies
$$\frac{\chi(F \otimes {\Cal M}_i
\otimes A^n)}{rk(F)} > \frac{\chi(E \otimes {\Cal M}_i \otimes
A^n)}{rk(E)} \text{\ for\ }n >> 0,$$
contradicting the ${\Cal M}_i$-twisted $A$-Gieseker-semistability of
$E$.

Thus $E \otimes {\Cal L}_i \in S(r,{c_1}_{{\Cal L}_i},{c_2}_{{\Cal
L}_i},A)_{H'}$.

On the other hand, suppose $E \otimes {\Cal L}_i \in S(r,{c_1}_{{\Cal
L}_i},{c_2}_{{\Cal
L}_i},A)_{H'}$.  Take a subsheaf $F \otimes {\Cal L}_i \subset E
\otimes {\Cal L}_i$.  Since $E \otimes {\Cal L}_i$ is
$A$-Gieseker-semistable, we have
$$\frac{\chi(F \otimes {\Cal L}_i \otimes A^n)}{rk(F)} \leq
\frac{\chi(E \otimes {\Cal L}_i \otimes A^n)}{rk(E)} \text{\ for\
}n >> 0.$$
Therefore,
$$\frac{c_1(F)}{rk(F)} \cdot A \leq \frac{c_1(E)}{rk(E)} \cdot A$$
and if the equality holds, then
$$\align
&\{\frac{c_1(F)}{rk(F)} - \frac{c_1(E)}{rk(E)}\} \cdot {\Cal L}_i\\
&+ \frac{1}{2}\frac{c_1(F)^2 - 2c_2(F) - c_1(F) \cdot K_X}{rk(F)}\\
&- \frac{1}{2}\frac{c_1(E)^2 - 2c_2(E) - c_1(E) \cdot K_X}{rk(E)}
\leq 0.\\
\endalign$$
If
$$\frac{c_1(F)}{rk(F)} \cdot A < \frac{c_1(E)}{rk(E)} \cdot A,$$
then clearly
$$\frac{\chi(F \otimes {\Cal M}_i \otimes A^n)}{rk(F)} <
\frac{\chi(E \otimes {\Cal M}_i \otimes A^n)}{rk(E)} \text{\ for\
}n >> 0.$$
If
$$\frac{c_1(F)}{rk(F)} \cdot A = \frac{c_1(E)}{rk(E)} \cdot A$$
and
$$\align
&\{\frac{c_1(F)}{rk(F)} - \frac{c_1(E)}{rk(E)}\} \cdot {\Cal L}_i\\
&+ \frac{1}{2}\frac{c_1(F)^2 - 2c_2(F) - c_1(F) \cdot K_X}{rk(F)}\\
&- \frac{1}{2}\frac{c_1(E)^2 - 2c_2(E) - c_1(E) \cdot K_X}{rk(E)}
< 0,\\
\endalign$$
then
$$\align
&\{\frac{c_1(F)}{rk(F)} - \frac{c_1(E)}{rk(E)}\} \cdot {\Cal M}_i\\
&+ \frac{1}{2}\frac{c_1(F)^2 - 2c_2(F) - c_1(F) \cdot K_X}{rk(F)}\\
&- \frac{1}{2}\frac{c_1(E)^2 - 2c_2(E) - c_1(E) \cdot K_X}{rk(E)}
< 0\\
\endalign$$
and thus
$$\frac{\chi(F \otimes {\Cal M}_i \otimes A^n)}{rk(F)} <
\frac{\chi(E \otimes {\Cal M}_i \otimes A^n)}{rk(E)} \text{\ for\
}n >> 0.$$
If
$$\frac{c_1(F)}{rk(F)} \cdot A = \frac{c_1(E)}{rk(E)} \cdot A$$
and
$$\align
&\{\frac{c_1(F)}{rk(F)} - \frac{c_1(E)}{rk(E)}\} \cdot {\Cal L}_i\\
&+ \frac{1}{2}\frac{c_1(F)^2 - 2c_2(F) - c_1(F) \cdot K_X}{rk(F)}\\
&- \frac{1}{2}\frac{c_1(E)^2 - 2c_2(E) - c_1(E) \cdot K_X}{rk(E)}
= 0,\\
\endalign$$
i.e., if
$$p(F \otimes {\Cal L}_i,A,a) = p(E \otimes {\Cal L}_i,A,a),$$
then
we have
$$\{\frac{c_1(F)}{rk(F)} - \frac{c_1(E)}{rk(E)}\} \cdot H' \geq 0.$$
If
$$\{\frac{c_1(F)}{rk(F)} - \frac{c_1(E)}{rk(E)}\} \cdot H' > 0,$$
then
$$\align
&\{\frac{c_1(F)}{rk(F)} - \frac{c_1(E)}{rk(E)}\} \cdot {\Cal M}_i\\
&+ \frac{1}{2}\frac{c_1(F)^2 - 2c_2(F) - c_1(F) \cdot K_X}{rk(F)}\\
&- \frac{1}{2}\frac{c_1(E)^2 - 2c_2(E) - c_1(E) \cdot K_X}{rk(E)}
< 0\\
\endalign$$
and hence
$$\frac{\chi(F \otimes {\Cal M}_i \otimes A^n)}{rk(F)} <
\frac{\chi(E \otimes {\Cal M}_i \otimes A^n)}{rk(E)} \text{\ for\
}n >> 0.$$
If
$$\{\frac{c_1(F)}{rk(F)} - \frac{c_1(E)}{rk(E)}\} \cdot H' = 0,$$
then
$$\frac{\chi(F \otimes {\Cal M}_i \otimes A^a)}{rk(F)} =
\frac{\chi(E \otimes {\Cal M}_i \otimes A^a)}{rk(E)} \text{\ for\
all\ }a.$$
Thus $E \in S((r,c_1,c_2) \otimes {\Cal M}_i,A)$.

A similar argument works for the other cases inside of the (resp.).
This proves (i).

\vskip.1in

First note that the existence of a good categorical quotient $M(r,{c_1}_{{\Cal
L}_i},{c_2}_{{\Cal
L}_i},A)_{H'}$ follows from the standard GIT [MFK94].  Since the
locus $Q$ is a closed subscheme of the projective Quot scheme
$Quot({\Cal O}_X^{\oplus l}/\chi_{r,{c_1}_{{\Cal L}_i}',{c_2}_{{\Cal
L}_i}'})$, the GIT quotient $M(r,{c_1}_{{\Cal L}_i},{c_2}_{{\Cal
L}_i},A)_{H'}$ is also a projective scheme.  Also note that
set-theoretically we can identify $S((r,c_1,c_2) \otimes {\Cal M}_i,
A)$ with $S(r,{c_1}_{{\Cal L}_i},{c_2}_{{\Cal L}_i},A)_{H'}$ by (i).
Now we verify that for ${\Cal M}_i$-twisted $A$-Gieseker-semistable
sheaves $E$ and $E'$, the closures of the orbits corresponding to $E
\otimes {\Cal L}_i \otimes A^a$ and $E' \otimes {\Cal L}_i \otimes
A^a$ in $Q^{ss}$ intersect if and only if $gr(E) \cong gr(E')$, where
$Q^{ss}$ is the locus of semistable points in $Q$ (See \S 2.) with
respect to the linearization induced from the Pl\"ucker embedding by
taking a high multiple of $H'$ and where $gr(E)$ is the direct sum of
the quotients of the Harder-Narasimhan filtration of $E$ with respect
to ${\Cal M}_i$-twisted $A$-Gieseker-stability.  Given an extension
$$0 \rightarrow {\Cal F} \rightarrow {\Cal E} \rightarrow {\Cal G}
\rightarrow 0$$ we can find a family of extensions ${\Cal E}_t$ of
${\Cal G}$ by ${\Cal F}$, parametrized by $t \in A_1^k$, such that
for each $t \neq 0$ the extension is the given one, and for $t = 0$
the extension is trivial.  Now applying this repeatedly to the
Harder-Narasimhan filtration we see that the orbit corresponding to
$gr(E) \otimes {\Cal L}_i \otimes A^a$ is in the closure of the orbit
corresponding to $E \otimes {\Cal L}_i \otimes A^a$.  Also note that
if $E$ is ${\Cal M}_i$-twisted $A$-Gieseker-semistable, then so is
$gr(E)$ and thus the orbit corresponding to $gr(E) \otimes {\Cal L}_i
\otimes A^a$ is in $Q^{ss}$.  So if $gr(E) = gr(E')$ then the
closures of the orbits of $E \otimes {\Cal L}_i \otimes A^a$ and $E'
\otimes {\Cal L}_i \otimes A^a$ intersect.

Conversely we see that the orbit corresponding to $gr(E) \otimes {\Cal L}_i
\otimes A^a$ is
closed as follows.  Suppose $R$ is a discrete valuation ring over
$k$ with the field of fractions $K$ and residue field $k$.  Suppose
${\Cal E}$ is a family of ${\Cal M}_i$-twisted
$A$-Gieseker-semistable sheaves on $X$ over $\text{Spec} R$ with
${\Cal E}_K \cong E \otimes K$, $E$ over $k$ with $gr(E) = E$.
If $E_i$ is a ${\Cal M}_i$-twisted $A$-Gieseker-stable component
of $E$, then by semicontinuity, then there are at least as many
maps from $E_i$ to ${\Cal E}_0$ as to $E$.  Since ${\Cal E}_0$ is
${\Cal M}_i$-twisted $A$-Gieseker-semistable, this implies that
${\Cal E}_0$ is a direct sum of copies of $E_i$ with the same
multiplicities as $E$, so ${\Cal E}_0 \cong E$.

The coarse representability for $M((r,c_1,c_2) \otimes
{\Cal M}_i,A)$ follows from the universal property of the Quot
scheme and the standard GIT.  Other cases follow similarly.  This
proves (ii).

(iii) is immediate from the Key GIT Lemma of \S 2 and (i) (ii)
above.  We remark that the morphism $\psi_i$ (resp. $\psi_i^+$) is
induced from the universal property of the moduli spaces
$M((r,c_1,c_2) \otimes {\Cal M}_i,A)$ (resp. $M((r,c_1,c_2) \otimes {\Cal
M}_{i+1},A)$) and $M((r,c_1,c_2) \otimes {\Cal L}_i,A)$, and that the
point in $M((r,c_1,c_2) \otimes {\Cal M}_i,A)$ (resp. $M((r,c_1,c_2)
\otimes {\Cal M}_{i+1},A)$) corresponding to the Seshadri equivalence
class of $E$ with respect to the ${\Cal M}_i$-twisted (resp.
${\Cal M}_{i+1}$-twisted) $A$-Gieseker-stability is mapped under
$\psi_i$ (resp. $\psi_i^+$) to the point corresponding to the Seshadri
equivalence class of $E$ with respect to the ${\Cal L}_i$-twisted
$A$-Gieseker-stability.

\vskip.1in

\S 5. Construction of Flip: Rational Case

\vskip.1in

In this section, we construct the flip
$$
\CD
M((r,c_1,c_2) \otimes {\Cal M}_i,A) @.@.@.@.M((r,c_1,c_2)
\otimes {\Cal M}_{i+1},A)\\
@.\searrow @.@.\swarrow @.\\
@.@.M((r,c_1,c_2) \otimes {\Cal L}_i,A)@.@.\\
\endCD
$$
in the more delicate case when a representative ${\Cal L}_i \in L_i$ is only a
RATIONAL line bundle ${\Cal L}_i \in \text{Pic}(X) \otimes {\Bbb Q}$.  Our
strategy is to use Kawamata's technique of finding some nice Galois cover
(which was originally used to prove the celebrated Kawamata-Viehweg
vanishing Theorem) $\phi:Y \rightarrow X$ so that the ${\Cal
L}_i$-twisted $A$-Gieseker-semistable sheaves on $X$ are pulled back
to the ${\Cal L}_i^Y$-twisted $\phi^*A$-Gieseker-semistable sheaves on
$Y$ , where ${\Cal L}_i^Y$ is now an INTEGRAL line
bundle on $Y$.  (Note that ${\Cal L}_i^Y$ is not exactly the pull
back $\phi^*{\Cal L}_i$ but ${\Cal L}_i^Y = \phi^*{\Cal L}_i +
\frac{1}{2}R$ as below where $R$ is the ramification divisor for
$\phi$.)  Then by looking at the appropriate locus of the Quot scheme
associated to $Y$ and taking the action of the Galois group into
consideration, we reduce the construction of the flip to the
Mumford-Thaddeus principle given by the Key GIT Lemma of \S 2 on $Y$.

\proclaim{Proposition 5.1}  Let ${\Cal L}_i =
\frac{p}{q}\Lambda_i$ where $\Lambda_i$ is a very ample line bundle
on $X$.  Take a positive integer $m \in {\Bbb N}$ such that $q$
divides $m$ and that $mA - \Lambda_i$ is very ample.  There exists a
Galois cover $\phi:Y \rightarrow X$ from a nonsingular projective
surface s.t.

\ \ (i) the Galois group $G \cong ({\Bbb Z}/m)^4$ and thus $\phi$ is
a Kummer extension,

\ \ (ii) ${\Cal L}_i^Y = \phi^*{\Cal L}_i + \frac{1}{2}R$ is
represented by an integral line bundle which is naturally a
$G$-sheaf in the sense of Mumford (See [Mumford70,P.69].), where $R$ is the
ramification divisor $K_Y = \phi^*K_X +
R$, and

\ \ (iii) a coherent sheaf on $X$ is torsion free and ${\Cal
L}_i$-twisted $A$-Gieseker-semistable if and only if $\phi^*E$ is
torsion free and ${\Cal L}_i^Y$-twisted $\phi^*A$-Gieseker-semistable.
\endproclaim

\demo{Proof of Proposition 5.1}\enddemo

First we construct the Galois cover following [KMM87,\S 1-1].

Take smooth irreducible members
$$\Gamma_1, \Gamma_2 \in
|\Lambda_i|$$
and
$$\align
H_1^{(1)}, H_2^{(1)} &\in |mA - \Gamma_1|\\
H_1^{(2)}, H_2^{(2)} &\in |mA - \Gamma_2|\\
\endalign$$
such that
$$\Gamma_1 \cup \Gamma_2 \cup \cup_{k,i}H_k^{(i)}$$
has only simple normal crossings.  Now let $X = \cup U_{\alpha}$ be
an affine open cover of $X$ with the transition functions of $A$
$$\{a_{\alpha\beta};
a_{\alpha\beta} \in H^0(U_{\alpha} \cap U_{\beta},{\Cal O}_X^*)\}$$
and local sections
$$\{\varphi_{k\alpha}^{(i)};
\varphi_{k\alpha}^{(i)} \in H^0(U_{\alpha},{\Cal O}_X)\}$$
such that
$$H_k^{(i)} + \Gamma_i|_{U_{\alpha}} =
\text{div}(\varphi_{k\alpha}^{(i)}) \text{\ on\ }U_{\alpha}$$
and that
$$\varphi_{k\alpha}^{(i)} = a_{\alpha\beta}^m \cdot
\varphi_{k\beta}^{(i)}.$$
We take the normalization of $X$ in
$\text{Rat}(X)[\cup_{k,i}(\varphi_{k\alpha}^{(i)}
)^{\frac{1}{m}}]$ (for some $\alpha$) as $Y$.
  (Note that $\text{Rat}(X)[\cup_{k,i}(\varphi_{k\alpha}^{(i)}
)^{\frac{1}{m}}] = \text{Rat}(X)[\cup_{k,i}(\varphi_{k\beta}^{(i)}
)^{\frac{1}{m}}]$ for any $\alpha, \beta$.)

[KMM87,Theorem 1-1-1 $\&$ Lemma 1-1-2] and the construction imply that
$\phi:Y \rightarrow X$ is a Galois cover from a nonsingular surface
$Y$ with the Galois group $G \cong ({\Bbb Z}/m)^4$ such that
$$\align
\phi^*(\Gamma_i) &= m\{\phi^*(\Gamma_i)_{red}\}\\
\phi^*(H_k^{(i)}) &= m\{\phi^*(H_k^{(i)})_{red}\}\\
\endalign$$
and that
$$K_Y = \phi^*K_X + R$$
where
$$\align
R &= (m - 1)\{\phi^*(\Gamma_1)_{red} + \phi^*(H_1^{(1)})_{red}
 + \phi^*(H_2^{(1)})_{red}\\
&+ \phi^*(\Gamma_2)_{red} +
\phi^*(H_1^{(2)})_{red} + \phi^*(H_2^{(2)})_{red}
\}\\
&= \frac{m - 1}{m}\{2\phi^*(2mA - \Gamma_1)\}.\\
\endalign$$
Therefore,
$$\align
{\Cal L}_i^Y &= \phi^*{\Cal L}_i + \frac{1}{2}R\\
&= \frac{p}{q}m\{(\phi^*\Gamma_1)_{red}\} + (m - 1)\phi^*(2mA) +
(m - 1)\{(\phi^*\Gamma_1)_{red}\}\\
\endalign$$
is integral since $q$ divides $m$.  Moreover, since
$(\phi^*\Gamma_1)_{red}$ is $G$-invariant, the line
bundle ${\Cal L}_i^Y$ can be given naturally the structure of a
$G$-sheaf in the sense of Mumford by embedding it in the constant
sheaf $\text{Rat}(X)$.  This proves (i) and (ii).

In order to prove (iii), first consider the exact sequence for a
coherent torsion free sheaf $E$ on $X$
$$0 \rightarrow E \rightarrow E^{**} \rightarrow Coker \rightarrow
0,$$
which gives rise to another exact sequence on $Y$
$$0 \rightarrow \phi^*E \rightarrow \phi^*E^{**}
\rightarrow \phi^*Coker \rightarrow 0,$$
since $\phi$ is faithfully flat.  But since $E^{**}$ is locally
free and thus so is $\phi^*E^{**}$, we conclude $\phi^*E$ is
torsion free.  The converse is also immediate.

Now the computation of the difference
$$\align
&\frac{\chi(\phi^*F \otimes {\Cal L}_i^Y \otimes
\phi^*A^n)}{rk(\phi^*F)} - \frac{\chi(\phi^*E \otimes {\Cal L}_i^Y
\otimes \phi^*A^n)}{rk(\phi^*E)}\\
&= \text{deg}\phi\{\frac{\chi(F \otimes {\Cal L}_i
\otimes A^n)}{rk(F)} - \frac{\chi(E \otimes {\Cal L}_i
\otimes A^n)}{rk(E)}\}\\
\endalign$$
shows that if $\phi^*E$ is ${\Cal L}_i^Y$-twisted
$\phi^*A$-Gieseker-semistable, then $E$ is ${\Cal L}_i^Y$-twisted
$A$-Gieseker-semistable.  To verify the converse we only have to
show that if $F_Y \subset \phi^*E$ is the first piece of the
Harder-Narasimhan filtration of $\phi^*E$ with respect to the
${\Cal L}_i^Y$-twisted $\phi^*A$-Gieseker-semistability, then $F_Y =
\phi^*F$ for some subsheaf $F \subset E$.

First remark that $\phi^*E$ is a $G$-sheaf and that the maximality
of $F_Y$ implies $F_Y$ is $G$-invariant and thus $F_Y$ is a
$G$-subsheaf of $\phi^*E$.  We will actually prove $F_Y =
\phi^*(\phi_*(F_Y)^G)$.

We claim
$$\phi^*(\phi_*(F_Y)^G) \subset F_Y \subset
\phi^*(\{\phi_*(F_Y)^G\}^s)$$
where $\phi_*(F_Y)^G \subset \{\phi_*(F_Y)^G\}^s \subset E$ is the
saturation of $\phi_*(F_Y)^G$ in $E$, and we will prove this by a
local analysis of the Galois cover $\phi:Y \rightarrow X$.

Take $x \in X$.  Set $R = {\Cal O}_{X,x}$.

We deal with the case where
$$\align
x &\in \Gamma_1, H_1^{(2)}\\
x &\notin \Gamma_2, H_2^{(1)}, H_1^{(1)}, H_2^{(2)}.\\
\endalign$$
(Other cases can be treated similarly.)

\vskip1.5in

{}From the construction of the Galois cover, for a set of local equations
$$\align
\varphi_1^{(1)} &\text{\ of\ }\Gamma_1 + H_1^{(1)}\\
\varphi_2^{(1)} &\text{\ of\ }\Gamma_1 + H_2^{(1)}\\
\varphi_1^{(2)} &\text{\ of\ }\Gamma_2 + H_1^{(2)}\\
\varphi_2^{(2)} &\text{\ of\ }\Gamma_2 + H_2^{(2)},\\
\endalign$$
we have
$$(\phi_*{\Cal O}_Y)_x = \oplus_{i,j,k,l}R \cdot
((\varphi_1^{(1)})^{\frac{1}{m}})^i
((\frac{\varphi_1^{(1)}}{\varphi_2^{(1)}})^{\frac{1}{m}})^j
((\varphi_1^{(2)})^{\frac{1}{m}})^k((\varphi_2^{(2)})^{\frac{1}{m}})^l,$$
This decomposition corresponds to the decomposition of $(\phi_*{\Cal O}_Y)_x$
into the eigen-\linebreak
spaces under the action of $G$.  (cf. [KMM87,Theorem
1-1-1 and Lemma 1-1-2].)  Now let $E_x = M$.  Then
$$\{\phi_*(\phi^*E)\}_x = \oplus_{i,j,k,l}M \otimes R
\cdot ((\varphi_1^{(1)})^{\frac{1}{m}})^i((\frac{\varphi_1^{(1)}}
{\varphi_2^{(1)}})^{\frac{1}{m}})^j
((\varphi_1^{(2)})^{\frac{1}{m}})^k
((\varphi_2^{(2)})^{\frac{1}{m}})^l.$$
Let
$$\align
\{\phi_*(\phi^*E)\}_x \supset \phi_*(F_Y)_x &=
\oplus_{i,j,k,l}{F_Y}_{i,j,k,l}
\otimes R
\cdot ((\varphi_1^{(1)})^{\frac{1}{m}})^i \\
& \qquad\cdot ((\frac{\varphi_1^{(1)}}
{\varphi_2^{(1)}})^{\frac{1}{m}})^j
((\varphi_1^{(2)})^{\frac{1}{m}})^k
((\varphi_2^{(2)})^{\frac{1}{m}})^l.\\
\endalign$$
Note that
$$\{\phi_*(F_Y)^G\}_x = {F_Y}_{0,0,0,0}.$$
Since $\phi_*(F_Y)_x$ is a $(\phi_*{\Cal O}_Y)_x$-module, we have
$$\align
{F_Y}_{0,0,0,0} &\subset {F_Y}_{i,j,k,l}\\
(\varphi_1^{(1)})(\frac{\varphi_1^{(1)}}{\varphi_2^{(1)}})(\varphi_1^{(2)})
(\varphi_2^{(2)}){F_Y}_{i,j,k,l} &\subset {F_Y}_{0,0,0,0}.\\
\endalign$$
This implies $F_Y/\phi^*(\phi_*(F_Y)^G)$ is torsion.  As
$$\phi_*(F_Y)^G \subset \{\phi_*(F_Y)^G\}^s$$
and
$$\phi^*E/\phi^*(\{\phi_*(F_Y)^G\}^s) \cong
 \phi^*(E/\{\phi_*(F_Y)^G\}^s)$$
is torsion free, we obtain the claim
$$F_Y \subset \phi^*(\{\phi_*(F_Y)^G\}^s).$$
Therefore, we have
$$\align
rk(\phi_*(F_Y)^G) &= rk(\phi^*(\phi_*(F_Y)^G) = rk(F_Y) \\
&=
rk(\phi^*(\{\phi_*(F_Y)^G\}^s)) = rk(\{\phi_*(F_Y)^G\}^s),\\
\endalign$$
and in particular
$$\frac{\chi(F_Y \otimes {\Cal L}_i^Y \otimes \phi^*A^n)}{rk(F_Y)} \leq
\frac{\phi^*
(\{\phi_*(F_Y)^G\}^s) \otimes {\Cal L}_i^Y \otimes
\phi^*A^n)}{rk(\phi^*(\{\phi_*(F_Y)^G\}^s))} \text{\ for\ }n >> 0.$$
By the maximality of $F_Y$, we obtain
$$F_Y = \phi^*(\{\phi_*(F_Y)^G\}^s),$$
and thus
$$\phi^*(\phi_*(F_Y)^G) = F_Y = \phi^*(\{\phi_*(F_Y)^G\}^s).$$
This completes the proof of (iii) and Proposition 5.1.

\vskip.2in

Now we take the Quot scheme
$$Quot_Y = Quot({\Cal O}_Y^{\oplus l_Y}/\chi_{(r,\phi^*c_1,\phi^*c_2) \otimes
{\Cal L}_i^Y \otimes \phi^*A^a})$$
which parametrizes the quotient of ${\Cal O}_Y^{\oplus l_Y}$ whose Hilbert
polynomial is the same as $\phi^*E \otimes {\Cal L}_i^Y \otimes \phi^*A^a$
for $E \in S((r,c_1,c_2) \otimes {\Cal L}_i,A)$ with $l_Y =
\chi(\phi^*E \otimes {\Cal L}_i^Y \otimes \phi^*A^a)$.
  Note that we take $a$ appropriately large according to the construction in \S
2 applied to this situation.

There is a natural action of $G$ on $Quot_Y$, namely if
$$\align
q:{\Cal O}_Y^{\oplus l_Y} &\rightarrow E_Y \rightarrow 0\\
(t_1, t_2, \cdot\cdot\cdot, t_{l_Y}) &\rightarrow (s_1, s_2, \cdot\cdot\cdot
s_{l_Y})\\
\endalign$$
is a point sending the sections $(t_1, t_2, \cdot\cdot\cdot, t_{l_Y})$
corresponding to the direct summands to the sections $(s_1, s_2,
\cdot\cdot\cdot,
s_{l_Y})$ of $E_Y$, then for $g \in G$
$$\align
gq:{\Cal O}_Y^{\oplus l_Y} &\rightarrow g^*E_Y \rightarrow 0\\
(t_1', t_2', \cdot\cdot\cdot, t_{l_Y}') &\rightarrow (g^*s_1, g^*s_2,
\cdot\cdot\cdot g^*s_{l_Y})\\
\endalign$$
is the point sending the sections $(t_1', t_2', \cdot\cdot\cdot, t_{l_Y}')$
corresponding to the direct summands to the sections $(g^*s_1, g^*s_2,
\cdot\cdot\cdot, g^*s_{l_Y})$ of $g^*E_Y$.  Note that we identify the element
$g
\in G$ with the automorphism $g:Y \rightarrow Y$ and that according to this
definition of the action, we have the reversed order identity $(gh)q = h(gq)$ a
priori, but having $G$ being abelian we have the usual identity for the action
 $(gh)q = g(hq)$.

There is also a natural action of $GL(l_Y)$ on $Quot_Y$.  If
$$\align
q:{\Cal O}_Y^{\oplus l_Y} &\rightarrow E_Y \rightarrow 0\\
(t_1, t_2, \cdot\cdot\cdot, t_{l_Y}) &\rightarrow (s_1, s_2, \cdot\cdot\cdot
s_{l_Y})\\
\endalign$$
is a point sending the sections $(t_1, t_2, \cdot\cdot\cdot, t_{l_Y})$
corresponding to the direct summands to the sections $(s_1, s_2,
\cdot\cdot\cdot,
s_{l_Y})$ of $E_Y$, then for $M \in GL(l_Y)$
$$\align
Mq:{\Cal O}_Y^{l_Y} (\rightarrow {\Cal O}_Y^{l_Y}) &\rightarrow E_Y \rightarrow
0\\
(t_1', t_2', \cdot\cdot\cdot, t_{l_Y}') &\rightarrow (s_1, s_2,
\cdot\cdot\cdot, s_{l_Y})M\\
\endalign
$$
is the point sending the sections $(t_1', t_2', \cdot\cdot\cdot, t_{l_Y}')$
corresponding to the direct summands to the sections $(s_1, s_2,
\cdot\cdot\cdot, s_{l_Y})M$.

Suppose
$$E \in S((r,c_1,c_2) \otimes {\Cal L}_i,A).$$
Then by Proposition 5.1 (iii)
$$\phi^*E \in S((r,\phi^*c_1,\phi^*c_2) \otimes {\Cal L}_i^Y,\phi^*A).$$
Moreover, Proposition 5.1 (ii) implies
$$\phi^*E \otimes {\Cal L}_i^Y \otimes \phi^*A^a = \phi^*(E \otimes A^a)
\otimes {\Cal L}_i^Y$$
has a natural $G$-sheaf structure, and thus $G$ acts on
$$H^0(\phi^*E \otimes {\Cal L}_i^Y \otimes \phi^*A^a).$$
We have the decomposition into eigenspaces
$$H^0(\phi^*E \otimes {\Cal L}_i^Y \otimes \phi^*A^a) = V = \oplus
V_{i,j,k,l}$$
under the action of $G$, and taking the basis from the eigenspaces we obtain
the
 matrix representation
$g \leftrightarrow L_g$, where $L_g$ is block diagonal with blocks of size
 $\dim V_{i,j,k,l}$,
corresponding to the simultaneous diagonalization of the matrix representation
 of $G$, i.e., for the basis
$$s_1, s_2, \cdot\cdot\cdot, s_{l_Y} \in H^0(\phi^*E \otimes {\Cal L}_i^Y
\otimes \phi^*A^a)$$
chosen from the eigenspaces $V_{i,j,k,l}$ we have
$$(gs_1, gs_2, \cdot\cdot\cdot, gs_{l_Y}) = (s_1, s_2, \cdot\cdot\cdot,
s_{l_Y})L_g.$$
Note that the above action of $g$ can be expressed as $\phi_g \circ g^*$ where
$$\phi^*E \otimes {\Cal L}_i^Y \otimes \phi^*A^a
\overset{\phi_g}\to{\overset{\sim}\to\leftarrow} g^*(\phi^*E \otimes {\Cal
L}_i^Y \otimes \phi^*A^a)$$
is the isomorphism satisfying the cocycle condition in the definition of the
$G$-sheaf structure of $\phi^*E \otimes {\Cal L}_i^Y \otimes \phi^*A^a$
(cf.[Mumford70]).

Remark that the dimension of $V_{i,j,k,l}$ is independent of
$E \in S((r,c_1,c_2) \otimes {\Cal L}_i,A)$ and if we fix the order of the
$V_{i,j,k,l}$ then $L_g$ is also independent of
 $E \in S((r,c_1,c_2) \otimes {\Cal L}_i,A)$ (at least on a connected component
of the Quot scheme, and we construct the moduli space connected component by
component).

Now take a point
$$\align
q:{\Cal O}_Y^{\oplus l_Y} &\rightarrow \phi^*E \otimes {\Cal L}_i^Y \otimes
\phi^*A^a \rightarrow 0\\ (t_1, t_2, \cdot\cdot\cdot, t_{l_Y}) &\rightarrow
(s_1,
s_2, \cdot\cdot\cdot s_{l_Y})\\ \endalign$$
where the sections $(s_1, s_2, \cdot\cdot\cdot, s_{l_Y})$ are chosen as above,
then the isomorphism
 $\phi_g$ makes the following diagram commutative
$$\align
gq:{\Cal O}_Y^{\oplus l_Y} &\rightarrow g^*(\phi^*E \otimes {\Cal L}_i^Y
\otimes \phi^*A^a) \rightarrow 0\\
(t_1', t_2', \cdot\cdot\cdot, t_{l_Y}') &\rightarrow (s_1, s_2, \cdot\cdot\cdot
s_{l_Y})\\
&\\
&\\
L_gq:{\Cal O}_Y^{l_Y} (\rightarrow {\Cal O}_Y^{l_Y}) &\rightarrow \phi^*E
\otimes {\Cal L}_i^Y \otimes \phi^*A^a \rightarrow 0\\
(t_1', t_2', \cdot\cdot\cdot, t_{l_Y}') &\rightarrow (s_1, s_2,
\cdot\cdot\cdot, s_{l_Y})M\\
\endalign$$
and thus
$$gq = L_gq \in Quot_Y.$$

\proclaim{Lemma 5.2} Define
$${\Cal D} = \{q \in Quot_Y;gq = L_gq \text{\ for\ }\forall g \in G\} \subset
Quot_Y.$$

\ \ (i) ${\Cal D}$ is naturally a closed subscheme of $Quot_Y$.

\ \ (ii) For
$$q:{\Cal O}_Y^{\oplus l_Y} \rightarrow E_Y \otimes {\Cal L}_i^Y \otimes
\phi^*A^a \rightarrow 0 \in {\Cal D}$$
by definition there exists an isomorphism $\phi_g$ which makes the following
diagram commutative
$$\align
gq:{\Cal O}_Y^{\oplus l_Y} &\rightarrow g^*(E_Y \otimes {\Cal L}_i^Y \otimes
\phi^*A^a) \rightarrow 0\\
&\\
&\\
L_gq:{\Cal O}_Y^{l_Y} (\rightarrow {\Cal O}_Y^{l_Y}) &\rightarrow E_Y \otimes
{\Cal L}_i^Y \otimes \phi^*A^a \rightarrow 0.\\
\endalign$$
Then $\{\phi_g;g \in G\}$ gives the structure of $G$-sheaf on $E_Y \otimes
{\Cal L}_i^Y \otimes \phi^*A^a$ and thus on $E_Y$.  Actually $\phi_g$ is the
restriction of the isomorphism $\Phi_g$
 defined over the universal quotient sheaf $(Univ)_{\Cal D}$ over ${\Cal D}$
and $\{\Phi_g;g \in G\}$ gives the structure of $G$-sheaf
 on $(Univ)_{\Cal D}$.
\endproclaim

\demo{Proof of Lemma 5.2}\enddemo

(i) ${\Cal D} = \cap (g \times L_g)^{-1}\Delta$ is naturally a closed subscheme
of $Quot_Y$, where $g \times L_g:Quot_Y \rightarrow Quot_Y \times Quot_Y$ is
the
product morphism and $\Delta \subset Quot_Y$ is the diagonal.

(ii) Let
$$\align
q:{\Cal O}_Y^{\oplus l_Y} &\rightarrow E_Y \otimes {\Cal L}_i^Y \otimes
\phi^*A^a \rightarrow 0 \in {\Cal D}\\
(t_1, t_2, \cdot\cdot\cdot, t_{l_Y}) &\rightarrow (s_1, s_2, \cdot\cdot\cdot
s_{l_Y})\\
\endalign$$
be a point in ${\Cal D}$ sending the sections $(t_1, t_2, \cdot\cdot\cdot,
t_{l_Y})$ corresponding to the direct summands
 to the sections $(s_1, s_2, \cdot\cdot\cdot, s_{l_Y})$ of $E_Y \otimes {\Cal
L}_i^Y \otimes \phi^*A^a$.  Then
$$\align
&\phi_{gh}((gh)^*s_1, (gh)^*s_2, \cdot\cdot\cdot, (gh)^*s_{l_Y})\\
&= (s_1, s_2, \cdot\cdot\cdot, s_{l_Y})L_{gh}\\
&= (s_1, s_2, \cdot\cdot\cdot, s_{l_Y})L_hL_g\\
&= \{\phi_h(h^*s_1, h^*s_2, \cdot\cdot\cdot, h^*s_{l_Y})\}L_g\\
&= \phi_g \circ g^*\{\phi_h(h^*s_1, h^*s_2, \cdot\cdot\cdot, h^*s_{l_Y})\}.\\
\endalign$$
Since $(s_1, s_2, \cdot\cdot\cdot, s_{l_Y})$ generate $E_Y \otimes {\Cal L}_i^Y
\otimes \phi^*A^a$, this shows that $\{\phi_g;g \in G\}$ satisfies the cocycle
condition and thus gives
 $E_Y \otimes {\Cal L}_i^Y \otimes \phi^*A^a$ the structure of a $G$-sheaf.
The assertion on $(Univ)_{\Cal D}$
 over ${\Cal D}$ can be checked similarly.

\vskip.1in

For $q:{\Cal O}_Y^{\oplus l_Y} \rightarrow (Univ)_q \rightarrow 0 \in {\Cal
D}$,
$(Univ)_q$ is a $G$-sheaf by Lemma 5.2 (ii) and hence $(Univ)_q \otimes {{\Cal
L}_i^Y}^{-1}$ is also given a
 natural $G$-sheaf structure.  There is a natural map
$$\phi^*\{\phi_*((Univ)_q \otimes {{\Cal L}_i^Y}^{-1})^G\} \rightarrow (Univ)_q
\otimes {{\Cal L}_i^Y}^{-1}.$$

We define the locus in ${\Cal D}$
$$\align
{\Cal P}^o = \{q \in {\Cal D};&\phi^*\{\phi_*((Univ)_q \otimes {{\Cal
L}_i^Y}^{-1})^G\} \rightarrow (Univ)_q \otimes {{\Cal L}_i^Y}^{-1}\\
&\text{\ is\ an\ isomorphism}\},\\
\endalign$$
which is easily seen to be open in ${\Cal D}$, and ${\Cal P}$ to be its (scheme
theoretic) closure in ${\Cal D}$.

Now let $Q_Y$ be the closure in $Quot_Y$ of points corresponding to torsion
free sheaves
with given chern classes as in \S 2.
  We define
$${\Cal V}{\Cal D} = {\Cal P} \cap Q_Y \subset {\Cal D},$$
which is $SL(\oplus GL(l_{i,j,k,l}))$-invariant, where $SL(\oplus
GL(l_{i,j,k,l}))$ is the subgroup of
 $SL(l_Y)$ formed by the elements which are block diagonal, with blocks
in $GL(l_{i,j,k,l})$.  Here,
$$l_{i,j,k,l} = \text{dim}_kV_{i,j,k,l}.$$

\proclaim{Theorem 5.3}
(i)
$$\align
{{\Cal V}{\Cal D}_{\phi^*H}}^{ss} &= {{Q_Y}_{\phi^*H}}^{ss} \cap {\Cal V}{\Cal
D}\\
(\text{resp.\ } {{\Cal V}{\Cal D}_{\phi^*A}}^{ss} &= {{Q_Y}_{\phi^*A}}^{ss}
\cap {\Cal V}{\Cal D}\\
{{\Cal V}{\Cal D}_{\phi^*H'}}^{ss} &= {{Q_Y}_{\phi^*H'}}^{ss} \cap {\Cal
V}{\Cal D})\\
\endalign
$$
where ${{Q_Y}_{\phi^*H}}^{ss}$ (resp. ${{Q_Y}_{\phi^*A}}^{ss}$,
${{Q_Y}_{\phi^*H'}}^{ss}$) is the locus of the semistable points of $Q_Y$ with
respect to the action of $SL(l_Y)$ and the linearization induced from the
Pl\"ucker embedding of $Grass(k^{\oplus l_Y} \otimes {\phi^*H}^m,R_m)$ (resp.
$Grass(k^{\oplus l_Y} \otimes {\phi^*A}^m,R_m)$, $Grass(k^{\oplus l_Y} \otimes
{\phi^*H'}^m,R_m)$) as in \S 2 (applied to $Y$ instead of $X$)
 and ${{\Cal V}{\Cal D}_{\phi^*H}}^{ss}$ (resp. ${{\Cal V}{\Cal
D}_{\phi^*A}}^{ss}$, ${{\Cal V}{\Cal D}_{\phi^*H'}}^{ss}$) is the locus of the
semistable points with respect to the action of $SL(\oplus GL(l_{i,j,k,l}))$
and
the same linearization ($SL(\oplus GL(l_{i,j,k,l})) \subset SL(l_Y)$).

(ii)
$$\align
M((r,c_1,c_2) \otimes {\Cal M}_i,A) &\cong {{\Cal V}{\Cal
D}_{\phi^*H'}}^{ss}//SL(\oplus GL(l_{i,j,k,l}))\\
M((r,c_1,c_2) \otimes {\Cal L}_i,A) &\cong {{\Cal V}{\Cal
D}_{\phi^*A}}^{ss}//SL(\oplus GL(l_{i,j,k,l}))\\
M((r,c_1,c_2) \otimes {\Cal M}_{i+1},A) &\cong {{\Cal V}{\Cal
D}_{\phi^*H}}^{ss}//SL(\oplus GL(l_{i,j,k,l})).\\
\endalign$$

(iii) The diagram of morphisms
$$
\CD
{{\Cal V}{\Cal D}_{\phi^*H'}}^{ss}//SL @.@.@.@.{{\Cal
V}{\Cal D}_{\phi^*H}}^{ss}//SL\\
@.\searrow @.@.\swarrow @.\\
@.@.{{\Cal V}{\Cal D}_{\phi^*A}}^{ss}//SL @.@.\\
\endCD
$$
where $SL=SL(\oplus GL(l_{i,j,k,l})$,
gives the desired flip
$$
\CD
M((r,c_1,c_2) \otimes {\Cal M}_i,A) @.@.@.@.M((r,c_1,c_2)
\otimes {\Cal M}_{i+1},A)\\
@.\psi_i \searrow @.@.\swarrow \psi_i^+@.\\
@.@.M((r,c_1,c_2) \otimes {\Cal L}_i,A).@.@.\\
\endCD
$$
\endproclaim

\demo{Proof of Theorem 5.3}\enddemo

(i) First note that
$${{\Cal V}{\Cal D}_{\phi^*H}}^{ss} \supset {{Q_Y}_{\phi^*H}}^{ss} \cap {\Cal
V}{\Cal D}$$
follows from the definition of semistable points, since $SL(l_Y)$-invariant
sections are automatically $SL(\oplus GL(l_{i,j,k,l}))$-invariant.  We show the
opposite inclusion as follows.

Let
$$V = \oplus V_{i,j,k,l}$$
be the decomposition of the vector space $V$ into the direct summands
$V_{i,j,k,l}$ of $\text{dim}_kV_{i,j,k,l} = l_{i,j,k,l}$,
$$v_{i,j,k,l}:V \rightarrow V_{i,j,k,l}$$
the projection onto the $\{i,j,k,l\}$-th factor.

The Hilbert-Mumford numerical criterion for semistability for the Grassmannian
$Grass(V \otimes W,R)$ of the $R$-dimensional quotients of the vector space
 $V \otimes W$ under the action of $SL(\oplus GL(l_{i,j,k,l}))$ with the
linearization induced
 by the Pl\"ucker embedding now reads:

\proclaim{Lemma 2.1'} If a point $p:V \otimes W \rightarrow U \rightarrow 0$ in
$Grass(V \otimes W,R)$ is
semistable for the action of $SL(\oplus GL(l_{i,j,k,l}))$ and the linearization
induced by the Pl\"ucker embedding, then for all nonzero subspaces $L \subset
V$
which preserves the decomposition, i.e., $L = \oplus v_i(L)$ we have
 $p(L \otimes W) \neq 0$ and
$$\frac{\text{dim}L}{\text{dim}\ p(L \otimes W)} \leq
\frac{\text{dim}V}{\text{dim}U}.$$  \endproclaim

Using this Lemma 2.1' we can prove the following slightly modified version of
Lemma 2.3 in a similar way
presented in \S 2.

\proclaim{Lemma 2.3'} For a fixed $a \in {\Bbb N}$ ($\geq a_0$) and an ample
line bundle $\phi^*H$ on $Y$, there exists
 $M_{\phi^*H} \in {\Bbb N}$ s.t. for all $m \geq M_{\phi^*H}$, the following
holds: If a point
$$q:{\Cal O}_Y^{\oplus l_Y} \rightarrow E_Y \otimes {\Cal L}_i^Y \otimes
\phi^*A^a \in {\Cal V}{\Cal D}$$
is semistable with respect to the action of $SL(\oplus GL(l_{i,j,k,l}))$ and
the linearization induced from the Pl\"ucker embedding
of $Grass(k^{\oplus l_Y} \otimes_k H^0(Y,\phi^*H^m,R_m))$, then the natural
homomorphism
$$k^{\oplus l_Y} = H^0({\Cal O}_Y^{\oplus l_Y}) \rightarrow H^0(E_Y
\otimes {\Cal L}_i^Y \otimes \phi^*A^a)$$
is injective, and for any nonzero $G$-sheaf quotient
$$E_Y \otimes {\Cal L}_i^Y \otimes \phi^*A^a \rightarrow G_Y \otimes {\Cal
L}_i^Y
\otimes \phi^*A^a \rightarrow 0$$
(the surjection is compatible with the $G$-sheaf structure of both $E_Y$ and
$G_Y$) with  $rk(G_Y) \neq 0$, we have
$$\frac{h^0(G_Y \otimes {\Cal L}_i^Y
\otimes \phi^*A^a)}{rk(G_Y)} \geq \frac{\chi(E_Y \otimes {\Cal L}_i^Y \otimes
\phi^*A^a)}{rk(E_Y)}.$$
Moreover, suppose that ${\Cal G}$ over $T \times Y$ is
a bounded family of coherent sheaves on $Y$ (independent of $a$).  Then there
 exists $a_2 \in {\Bbb N} (\geq a_0)$ such that for all $a \geq a_2$ and an
ample line bundle $\phi^*H$ on $Y$ the following holds: If a point
$$q:{\Cal O}_Y^{\oplus l_Y} \rightarrow E_Y \otimes {\Cal L}_i^Y \otimes
\phi^*A^a \rightarrow 0 \in {\Cal V}{\Cal D}$$
is semistable with respect to the action of $SL(\oplus GL(l_{i,j,k,l}))$ and
the linearization induced from the Pl\"ucker embedding of
 $Grass(k^{\oplus l_Y} \otimes_k H^0(Y,\phi^*H^m),R_m)$ and if the natural
homomorphism
$$k^{\oplus l_Y} = H^0({\Cal O}_Y^{\oplus l_Y}) \rightarrow H^0(E_Y \otimes
{\Cal L}_i^Y \otimes \phi^*A^a)$$
is an isomorphism, then for any nonzero $G$-sheaf quotient
$$E_Y \otimes {\Cal L}_i^Y \otimes \phi^*A^a \rightarrow G_Y \otimes {\Cal
L}_i^Y \otimes \phi^*A^a \rightarrow 0$$
(the surjection is compatible with the $G$-sheaf structure of both $E_Y$ and
$G_Y$) with $G_Y \cong {\Cal G}_t$ for some $t \in T$,
we have
$$\frac{h^0(G_Y \otimes {\Cal L}_i^Y \otimes \phi^*A^a)}{\chi(G_Y \otimes {\Cal
L}_i^Y \otimes \phi^*A^a \otimes \phi^*H^m)} \geq \frac{\chi(E_Y \otimes {\Cal
L}_i^Y \otimes \phi^*A^a)}{\chi(E_Y \otimes {\Cal L}_i^Y \otimes \phi^*A^a
\otimes \phi^*H^m)}.$$
\endproclaim

\vskip.1in

Now suppose
$$q:{\Cal O}_Y^{\oplus l_Y} \rightarrow E_Y \otimes {\Cal L}_i^Y \otimes
\phi^*A^a \in
{{\Cal V}{\Cal D}_{\phi^*H}}^{ss},$$
i.e., $q \in {\Cal V}{\Cal D}$ is semistable with respect to the action of
$SL(\oplus GL(l_{i,j,k,l}))$ and the same linearization induced by a multiple
of
$\phi^*H$ as for the action of $SL(l_Y)$.  If we assume that $E_Y$ has a
torsion
$T_Y$, then noting that $T_Y$ is $G$-invariant and thus $G$-subsheaf of $E_Y$
and using Lemma 2.3' we derive a contradiction in the same way as in the proof
of
``only if" part of (ii) in Key GIT Lemma of \S 2.  If we assume that $E_Y$ is
not ${\Cal L}_i^Y$-twisted $\phi^*A$-Gieseker-semistable, then noting that the
maximal destabilizing subsheaf $F_Y$  is $G$-invariant and thus a $G$-subsheaf
and that the quotient $E_Y/G_Y$ is a $G$-sheaf with the surjection $E_Y
\rightarrow E_Y/F_Y \rightarrow 0$  compatible with the $G$-sheaf structure,
again we derive a contradiction in the same way as in the proof of ``only if"
part of (ii) in Key GIT Lemma of \S 2.

Remark that we can also prove in a similar way to the proof of (i) in Theorem
4.1 (by retaking ${\Cal L}_i, {\Cal M}_i, A, H, H'$ we may assume they are all
on a line in $V(\Delta_s)$) that
$$E_Y \in S((r,\phi^*c_1,\phi^*c_2) \otimes {\Cal L}_i^Y,\phi^*A)_{\phi^*H}$$
if and only if
$$E_Y \in S((r,\phi^*c_1,\phi^*c_2) \otimes {\Cal M}_i^Y,\phi^*A).$$

Now suppose
$$E_Y \notin S((r,\phi^*c_1,\phi^*c_2) \otimes {\Cal
L}_i^Y,\phi^*A)_{\phi^*H}.$$
Then
$$E_Y \notin S((r,\phi^*c_1,\phi^*c_2) \otimes {\Cal M}_i^Y,\phi^*A).$$
Take $F_Y$ to be the maximal destabilizing subsheaf.  Then $F_Y$ is
$G$-invariant and thus a $G$-subsheaf.  (We may also assume that
 $F_Y$ has the same averaged Euler characteristics
$$\frac{\chi(F_Y \otimes {\Cal L}_i^Y \otimes \phi^*A^a)}{rk(F_Y)} =
\frac{\chi(E_Y \otimes {\Cal L}_i^Y \otimes \phi^*A^a)}{rk(E_Y)}$$
(as polynomials in $a$).  The rest of the argument goes without change as in
the last part of the proof of ``only if" part of (ii) in Key GIT Lemma of \S 2
to
 derive a contradiction.  Therefore, we conclude
$$E_Y \in S((r,\phi^*c_1,\phi^*c_2) \otimes {\Cal L}_i^Y,\phi^*A)_{\phi^*H},$$
and thus
$${{\Cal V}{\Cal D}_{\phi^*H}}^{ss} = {{Q_Y}_{\phi^*H}}^{ss} \cap {\Cal V}{\Cal
D}.$$
The other cases in (resp.) can be proved similarly.

\vskip.2in

(ii) First we prove the following lemma.

\proclaim{Lemma 5.4}
For
$$q:{\Cal O}_Y^{\oplus l_Y} \rightarrow (Univ)_q \rightarrow 0 \in {{\Cal
V}{\Cal D}_{\phi^*A}}^{ss},$$
(and thus automatically for
$$q:{\Cal O}_Y^{\oplus l_Y} \rightarrow (Univ)_q \rightarrow 0 \in {{\Cal
V}{\Cal D}_{\phi^*H}}^{ss}$$
and
$$q:{\Cal O}_Y^{\oplus l_Y} \rightarrow (Univ)_q \rightarrow 0 \in {{\Cal
V}{\Cal D}_{\phi^*H'}}^{ss})$$
the natural morphism
$$\phi^*\{\phi_*((Univ)_q \otimes {{\Cal L}_i^Y}^{-1})^G\} \rightarrow (Univ)_q
\otimes {{\Cal L}_i^Y}^{-1}$$
is an isomorphism.
\endproclaim

\demo{Proof of Lemma 5.4}\enddemo

Since by definition,
$${\Cal V}{\Cal D} = {\Cal P} \cap Q_Y \subset {\Cal D}$$
where ${\Cal P}$ is the closure of the locus
$$\align
{\Cal P}^o = \{q \in {\Cal D};&\phi^*\{\phi_*((Univ)_q \otimes {{\Cal
L}_i^Y}^{-1})^G\} \rightarrow (Univ)_q \otimes {{\Cal L}_i^Y}^{-1}\\
&\text{\ is\ an\ isomorphism}\},\\
\endalign$$
we only have to show that in ${{\Cal V}{\Cal D}_{\phi*A}}^{ss}$ the condition
$$\phi^*\{\phi_*((Univ)_q \otimes {{\Cal L}_i^Y}^{-1})^G\} \rightarrow (Univ)_q
\otimes {{\Cal L}_i^Y}^{-1}$$
being an isomorphism is a closed condition in parameter $q$.

Note first that for $q \in {{\Cal V}{\Cal D}_{\phi*A}}^{ss}$ since $(Univ)_q$
is torsion free, the natural morphism
$$\phi^*\{\phi_*((Univ)_q \otimes {{\Cal L}_i^Y}^{-1})^G\} \rightarrow (Univ)_q
\otimes {{\Cal L}_i^Y}^{-1}$$
is injective.  Suppose that
$$\phi^*\{\phi_*((Univ)_s \otimes {{\Cal L}_i^Y}^{-1})^G\} \hookrightarrow
(Univ)_s \otimes {{\Cal L}_i^Y}^{-1}$$
is an isomorphism for all $s \in S$ but possibly at $s_0 \in S$, where $S$ is
any one parameter subspace in ${{\Cal V}{\Cal D}_{\phi*A}}^{ss}$.

We have the decomposition
$$(Id_S \times \phi)_*((Univ)_S \otimes p_2^*({{\Cal L}_i^Y}^{-1})) =
\oplus_{i,j,k,l}{\Cal E}_{i,j,k,l}$$
into the eigenspaces under the action of $G$ and
$$(Id_S \times \phi)_*((Univ)_S \otimes p_2^*({{\Cal L}_i^Y}^{-1}))^G = {\Cal
E}_{0,0,0,0}.$$
Since $(Id_S \times \phi_*)_S((Univ)_S)$ is flat over $S$,
$(Id_S \times \phi)_*((Univ)_S \otimes p_2^*({{\Cal L}_i^Y}^{-1}))$ is also
flat
over $S$.  This in turn implies $(Id_S \times \phi)_*((Univ)_S \otimes
p_2^*({{\Cal L}_i^Y}^{-1}))^G = {\Cal E}_{0,0,0,0}$ is flat over $S$.
Therefore,
we conclude that
 $(Id_S \times \phi)^*((Id_S \times \phi)_*((Univ)_S \otimes p_2^*{{\Cal
L}_i^Y}^{-1})^G)$ is flat over $S$.
Now since
$$(Univ)_{s_0} \otimes {{\Cal L}_i^Y}^{-1}/\phi^*(\phi_*\{(Univ)_{s_0} \otimes
{{\Cal L}_i^Y}^{-1}\}^G)$$
is torsion (cf.Proposition 5.1 (iii)) and since
$$\align
\chi((Univ)_{s_0} \otimes {{\Cal L}_i^Y}^{-1}&  \otimes {\Cal L}_i^Y \otimes
\phi^*A^a) = \chi((Univ)_{s} \otimes {{\Cal L}_i^Y}^{-1} \otimes {\Cal L}_i^Y
\otimes \phi^*A^a)\\
&= \chi(\phi^*(\phi_*\{(Univ)_{s} \otimes {{\Cal L}_i^Y}^{-1}\}^G) \otimes
{\Cal L}_i^Y \otimes \phi^*A^a)\\
&= \chi((Id_S \times \phi)_*((Univ)_S \otimes {{\Cal L}_i^Y}^{-1})^G|_s \otimes
{\Cal L}_i^Y \otimes \phi^*A^a)\\
&= \chi((Id_S \times \phi)_*((Univ)_S \otimes {{\Cal L}_i^Y}^{-1})^G|_{s_0}
\otimes {\Cal L}_i^Y \otimes \phi^*A^a)\\
&= \chi(\phi^*\{\phi_*((Univ)_{s_0} \otimes {{\Cal L}_i^Y}^{-1})^G\} \otimes
{\Cal L}_i^Y \otimes \phi^*A^a),\\ \endalign$$
we conclude
$$\phi^*\{\phi_*((Univ)_{s_0} \otimes {{\Cal L}_i^Y}^{-1})^G\} = (Univ)_{s_0}
\otimes {{\Cal L}_i^Y}^{-1}.$$

\vskip.1in

Lemma 5.4 implies that for
$$q:{\Cal O}_Y^{\oplus l_Y} \rightarrow E_Y \otimes {\Cal L}_i^Y \otimes
\phi^*A^a \in {{\Cal V}{\Cal D}_{\phi^*A}}^{ss},$$
we have
$$E_Y = \phi^*(\phi_*(E_Y)^G),$$
and that moreover
$$\align
&(Id \times \phi)^*\{(Id \times \phi)_*{((Univ) \otimes p_2^*({\Cal L}_i^Y
\otimes \phi^*A^a)^{-1})_{{{\Cal V}{\Cal D}_{\phi^*A}}^{ss}\times X}}^G\}\\ &
\rightarrow ((Univ) \otimes p_2^*({\Cal L}_i^Y \otimes \phi^*A^a)^{-1})_{{{\Cal
V}{\Cal D}_{\phi^*A}}^{ss} \times X}\\ \endalign$$
is an isomorphism.

In order to prove (ii)
$$M((r,c_1,c_2) \otimes {\Cal M}_i,A) \cong {{\Cal V}{\Cal
D}_{\phi^*H}}^{ss}//SL(\oplus GL(l_{i,j,k,l})),$$   we first verify that for
$$q:{\Cal O}_Y^{\oplus l_Y} \rightarrow E_Y \otimes {\Cal L}_i^Y \otimes
\phi^*A^a \in {{\Cal V}{\Cal D}_{\phi^*H}}^{ss}$$
and
$$q':{\Cal O}_Y^{\oplus l_Y} \rightarrow {E_Y}' \otimes {\Cal L}_i^Y \otimes
\phi^*A^a \in {{\Cal V}{\Cal D}_{\phi^*H}}^{ss},$$
the closures of the orbits of $q$ and $q'$ in ${\Cal V}{\Cal D}_{\phi^*H}^{ss}$
intersect if and only if $gr(\phi_*(E_Y)^G) = gr(\phi_*({E_Y}')^G)$ where
$gr(\phi_*(E_Y)^G)$ (resp. $gr(\phi_*({E_Y}')^G)$) is the direct sum of the
quotients of the Harder-Narasimhan filtration of $\phi_*(E_Y)^G$ (resp.
$\phi_*({E_Y}')^G$) with respect to the  ${\Cal M}_i$-twisted
$A$-Gieseker-stability.  The proof goes without any change from that for (ii)
in
Theorem 4.1.  Secondly if we have a family ${\Cal E}$ of ${\Cal M}_i$-twisted
$A$-Gieseker-semistable sheaves on $S \times X$, then by taking $(Id_S \times
\phi)^*({\Cal E} \otimes p_2^*({\Cal L}_i^Y \otimes \phi^*A^a))$ over $S \times
Y$ and locally choosing a frame for ${p_1}_*\{(Id_S \times \phi)^*({\Cal E}
\otimes p_2^*({\Cal L}_i^Y \otimes \phi^*A^a))\}$ which corresponds to taking
the
eigenspace decomposition, we have a map $$S \rightarrow {{\Cal V}{\Cal
D}_{\phi^*H}}^{ss}$$ unique up to the action of $SL(\oplus GL(l_{i,j,k,l})$.
Therefore, it gives a map  $$S \rightarrow {{\Cal V}{\Cal
D}_{\phi^*H}}^{ss}//SL(\oplus GL(l_{i,j,k,l})).$$   By uniqueness the map is
defined globally
 on $S$.  ${{\Cal V}{\Cal D}_{\phi^*H}}^{ss}//SL(\oplus GL(l_{i,j,k,l}))$ is
universal for this property since it is a categorical quotient, and this proves
$$M((r,c_1,c_2) \otimes {\Cal M}_i,A) \cong {\Cal V}{\Cal
D}_{\phi^*H}^{ss}//SL(\oplus GL(l_{i,j,k,l})).$$
The other cases in (resp.) can be proved similarly.

(iii) is now immediate from (i) (ii) above.  We remark that the morphism
$\psi_i$ (resp. $\psi_i^+$) is
induced from the universal property of the moduli spaces
$M((r,c_1,c_2) \otimes {\Cal M}_i,A)$ (resp. $M((r,c_1,c_2) \otimes {\Cal
M}_{i+1},A)$) and $M((r,c_1,c_2) \otimes {\Cal L}_i,A)$, and that the point in
$M((r,c_1,c_2) \otimes {\Cal M}_i,A)$ (resp. $M((r,c_1,c_2) \otimes {\Cal
M}_{i+1},A)$) corresponding to the Seshadri equivalence class of $E$
with respect to the ${\Cal M}_i$-twisted (resp.
${\Cal M}_{i+1}$-twisted) $A$-Gieseker-stability is mapped under
$\psi_i$ (resp. $\psi_i^+$) to the point corresponding to the Seshadri
equivalence class of $E$ with respect to the ${\Cal L}_i$-twisted
$A$-Gieseker-stability.

This completes the proof of Theorem 5.3.

\vskip.1in

Combining Theorem 4.1 and Theorem 5.3 we have

\proclaim{Theorem 5.4} Let $M_0, L_0, L_1, M_1, \cdot\cdot\cdot, L_l, M_{l+1}$
be a
sequence of strata starting with $M_0$ containing $H^n$ (for $n >>
0$) and ending with $M_{l+1}$ containing ${H'}^{n'}$ (for $n' >>
0)$ such that
$$\overline{M_i} \cap \overline{M_{i+1}} = L_i$$
for $i = 0, 1, \cdot\cdot\cdot, l$.  (Note that these strata
actually exhaust all the strata in $V(\Delta_s) = \coprod L_i
\coprod M_j$.) Then
$$\align
M(r,c_1,c_2,H) &= M((r,c_1,c_2) \otimes {\Cal M}_0,A)\\
M(r,c_1,c_2,H') &= M((r,c_1,c_2) \otimes {\Cal M}_{l+1},A)\\
M(r,c_1,c_2,A) &= M((r,c_1,c_2) \otimes {\Cal M}_i,A) \text{\
or\ }M((r,c_1,c_2) \otimes {\Cal L}_i,A) \text{\ for\ some\
}i,\\
\endalign$$
and there is a sequence of Thaddeus-type
flips
$$
\CD
M((r,c_1,c_2) \otimes {\Cal M}_i,A) @.@.@.@.M((r,c_1,c_2)
\otimes {\Cal M}_{i+1},A)\\
@.\searrow @.@.\swarrow @.\\
@.@.M((r,c_1,c_2) \otimes {\Cal L}_i,A)@.@.\\
\endCD
$$
for $i = 1, 2, \cdot\cdot\cdot, l$, each of which is
a transformation constructed by the Key GIT Lemma of \S 2
and thus governed by the Mumford-Thaddeus principle.
\endproclaim

\proclaim{Remark 5.5}\endproclaim

The structure of each flip
$$
\CD
M((r,c_1,c_2) \otimes {\Cal M}_i,A) @.@.@.@.M((r,c_1,c_2)
\otimes {\Cal M}_{i+1},A)\\
@.\psi_i \searrow @.@.\swarrow \psi_i^+@.\\
@.@.M((r,c_1,c_2) \otimes {\Cal L}_i,A)@.@.\\
\endCD
$$
can be fairly explicitly described.  For example, if $r = 2$, then over a point
$$p \in M((r,c_1,c_2) \otimes {\Cal L}_i,A)$$
where $\psi_i$ is not isomorphic, $p$ corresponds to the Seshadri equivalence
class of the form
$$E = L \oplus L'$$
with $L$ and $L'$ being rank one torsion free sheaves s.t. $E$ is ${\Cal
L}_i$-twisted $A$-Gieseker-semistable, but $L$ destabilizes $E$ with respect to
the ${\Cal M}_{i+1}$-twisted $A$-Gieseker-\linebreak semistability
 and $L'$ destabilizes $E$ with respect to ${\Cal M}_{i}$-twisted
$A$-Gieseker-\linebreak semistability.

We have the morphisms which are set-theoretically bijective
$$\align
{\Bbb P}(\text{Ext}^1(L,L')) &\rightarrow \psi_i^{-1}(p)\\
{\Bbb P}(\text{Ext}^1(L',L)) &\rightarrow {\psi_i^+}^{-1}(p),\\
\endalign$$
i.e., the flip corresponds to flipping the factors of the extension class (at
least set-theoretically).  More detailed and scheme-theoretic analysis of the
flip will be published elsewhere.

\vskip.1in

Finally starting from the $d$-cells of maximal dimension and descending
inductively on the dimension $d$, i.e., going from
$\Delta_s$ to $W$ and repeating, and moreover letting $\Delta$ vary in
$Amp(X)_{\Bbb Q}$, we have the main theorem.

\proclaim{Theorem 5.6 = Main Theorem} The moduli space $M(r,c_1,c_2,H)$ goes
through a sequence of flips (and contraction morphisms
 $\&$ their inverses) in the category of moduli spaces $M((r,c_1,c_2) \otimes
{\Cal L},A)$ of ${\Cal L}$-twisted $A$-Gieseker-semistable sheaves
 for rational line bundles ${\Cal L} \in \text{Pic}(X) \otimes {\Bbb Q}$ and $A
\in Amp(X)_{\Bbb Q}$, all of which are governed by
 the Mumford-Thaddeus principle of GIT.
\endproclaim

Though we restricted ourselves to the subvector space $V(\Delta_s) \subset
N^1(X)_{\Bbb Q}$ in Proposition 3.5 and Theorem 5.4,
 it can be shown in the same manner that there exists a stratification of
$N^1(X)_{\Bbb Q}$ which describes the change of ${\Cal L}$-twisted
$A$-Gieseker-semistable sheaves
 when ${\Cal L}$ varies in $N^1(X)_{\Bbb Q}$, and the corresponding moduli
spaces are connected by sequences of flips governed
by Mumford-Thaddeus principle.

\proclaim{Theorem 5.7} The moduli spaces $M((r,c_1,c_2) \otimes {\Cal L},A)$ of
${\Cal L}$-twisted $A$-Gieseker-semistable sheaves
 exist for rational line bundles ${\Cal L} \in \text{Pic}(X) \otimes {\Bbb Q}$
and polarizations $A \in Amp(X)_{\Bbb Q}$ and are connected by sequences of
flips
 (and contraction morphisms $\&$ their inverses), all of which are governed by
Mumford-Thaddeus principle of GIT.
\endproclaim

\vskip.1in

Note added in proof:
While revising the paper after circulating the first version in July 1994, we
learned in October 1994 that [Ellingsrud-G\"ottsche94] obtained similar
results in the rank 2 case on particular (although from the point of view
of Donaldson theory the most interesting) types of surfaces, but with a
much finer
analysis of flips.  Though our results have no
restrictions on surfaces or rank, the price had to be paid in getting only a
general description.  The aforementioned paper uses the notion of
parabolic semistability, and the construction of moduli spaces depends
on [Maruyama-Yokogawa93][Yokogawa93].  Their semistability coincides with our
notion of rationally-twisted Gieseker semistability as we will
explain below.  We also direct the reader's  attention to the recent
papers [Friedman-Qin94] and [Yoshioka94].

We use the notation in \S 3 ,4 and 5.  Let $H \in
\Delta_s$
be an ample line bundle.  Then a sufficiently high multiple $H^n$ is
in
the stratum $M_0$ whereas $H^{-n}$ is in the stratum $M_{l+1}$.  We take an
effective Cartier divisor $D \in H^n$.  Parabolic semistability
with
parameter $a \in {\Bbb Q}$ is measured by the polynomial $$Par_a(E,n) = (1 -
a)\frac{\chi(E \otimes {\Cal O}_X \otimes A^n)}{rk(E)} + a\frac{\chi(E \otimes
{\Cal O}_X(D) \otimes A^n)}{rk(E)}$$ whereas our ${\Cal L}$-twisted
$A$-Gieseker-semistability for ${\Cal L} = (1 - a)H^{-n} + aH^n \in
\text{Pic}(X) \otimes {\Bbb Q}$ is measured by the polynomial
$$Twist_{\Cal L}(E,n) = \frac{\chi(E \otimes {\Cal L} \otimes
A^n)}{rk(E)}.$$
Now it is straightforward to see that for a subsheaf $F \subset E$
$$Par_a(E,n) - Par_a(F,n) = Twist_{\Cal L}(E,n) - Twist_{\Cal L}(F,n)$$
and thus both semistabilities coincide.  In any case, we observe the
change of semistability according to the stratification
$$V(\Delta_s) = \coprod L_i \coprod M_j$$
starting from the point $H^n$ to the point $H^{-n}$ along the line joining
them.

\vskip.1in

\S 6. References.

\vskip.1in

[Bogomolov79]; F.A. Bogomolov, Holomorphic tensor and vector bundles
on projective varieties, Math. U.S.S.R. Izvestiya 13 (1979), 499-555

[CKM88]; H. Clemens-J. Koll\'ar-S. Mori, Higher Dimensional Complex Geometry,
Ast\'erisque vol. 166,
 Soc. Math. France (1988)

[Dolgachev-Hu93]; I. Dolgachev-Y. Hu, Variation of Geometric Invariant Theory
Quotient, preprint (1993)

[Ellingsrud-G\"ottsche94]; G. Ellingsrud-L. G\"ottsche, Variation of moduli
spaces and Donaldson invariants under change of polarization, preprint (1994)

[Friedman-Qin94]; R. Friedman-Z. Qin, Flips of moduli spaces and transition
formulas for Donaldson polynomial invariants of rational surfaces, preprint
(1994)

[Gieseker77]; D. Gieseker, On the moduli of vector bundles on an algebraic
surface, Annals of Math. 106 (1977),
 45-60

[Grothendieck60]; A. Grothendieck, Techniques de construction et th\'eor\`emes
d'existence en g\'eom\'etrie
 alg\'ebrique IV; les sch\'emas de Hilbert, Sem. Bourbaki Expose 221 (1960)

[Hartshorne77]; R. Hartshorne, Algebraic Geometry, Graduate Texts in Math. 52,
 Springer-Verlag (1977)

[Hu-W.Li93]; Variation of the Gieseker and Uhlenbeck compactifications, pre-
\linebreak
print (1993)

[J.Li93]; J. Li, Algebraic geometric interpretation of Donaldson's polynomial
invariants, J. Diff. Geom. 37 (1993), 417-466

[KMM87]; Y. Kawamata-K. Matsuda-K. Matsuki, Introduction to the minimal model
problem, Algebraic Geometry, Sendai, Adv. Studies in Pure Math.,
 vol. 10 T. Oda ed., Kinokuniya-North-Holland (1987), 283-360

[Kov\'acs94]; S. Kov\'acs, The cone of curves of a K3 surface, preprint (1994)

[Maruyama-Yokogawa92] M. Maruyama-K. Yokogawa, Moduli of parabolic stable
sheaves, Math. Ann. 293 (1992), 77-99

[Mori88]; S. Mori, Flip theorem and the existence of minimal models for
3-folds, J. Amer. Math. Soc. 1 (1988), 117-253

[Mumford70]; D. Mumford, Abelian Varieties, Tata Studies in Math., Oxford Univ.
Press (1970)

[MFK94]; D. Mumford-J. Fogarty-F. Kirwan, Geometric Invariant Theory 3rd
edition, Ergebnisse der Mathmatik und ihrer Grentzbiete 34, Springer-Verlag
(1994)

[Simpson92]; C. Simpson, Moduli of representations of the fundamental groups of
a smooth projective variety I and II, preprint (1992)

[Thaddeus93]; M. Thaddeus, Stable pairs, linear systems and the Verlinde
formula, preprint (1993)

[Thaddeus94]; M. Thaddeus, Geometric invariant theory and flips, preprint
(1994)

[Yokogawa93]; K. Yokogawa, Compactification of moduli of parabolic sheaves
and moduli for parabolic Higgs sheaves, J. Math. Kyoto University 33 (1993),
451-504

[Yoshioka94]; K. Yoshioka, Chamber structure of polarizations and
the moduli of stable sheaves on a ruled surface, preprint (1994)

\enddocument